%
\magnification=\magstep1   
\input amstex
\UseAMSsymbols
\input pictex
\vsize=23truecm
\NoBlackBoxes
\parindent=18pt
  
   \font\rmk=cmr8    \font\itk=cmti8  \font\ttk=cmtt8

\font\gross=cmbx10 scaled\magstep1 

\def\mod{\operatorname{mod}}

\def\Hom{\operatorname{Hom}}

\def\Ext{\operatorname{Ext}}
\def\rad{\operatorname{rad}}

\def\arr#1#2{\arrow <1.5mm> [0.25,0.75] from #1 to #2}

   

 \def\Rahmenbio#1%
   {$$\vbox{\hrule\hbox%
                  {\vrule%
                       \hskip0.5cm%
                            \vbox{\vskip0.3cm\relax%
                               \hbox{$\displaystyle{#1}$}%
                                  \vskip0.3cm}%
                       \hskip0.5cm%
                  \vrule}%
           \hrule}$$}        

\def\Rahmen#1%
   {\centerline{\vbox{\hrule\hbox%
                  {\vrule%
                       \hskip0.5cm%
                            \vbox{\vskip0.3cm\relax%
                               \hbox{{#1}}%
                                  \vskip0.3cm}%
                       \hskip0.5cm%
                  \vrule}%
           \hrule}}}

\vglue3cm
\centerline{\gross Cluster-additive functions on stable translation quivers}
                    		    	          	                      \bigskip
\centerline{Claus Michael Ringel}     
						          \bigskip\bigskip
\plainfootnote{}
{\rmk 2010 \itk Mathematics Subject Classification. \rmk 
Primary:
 16G70. 
Secondary:
 13F60, 
 16G20, 
 16D90. 
}

{\narrower \rmk Abstract.
Additive functions on translation quivers have played an important role
in the representation theory of finite dimensional algebras, the most
prominent ones are the hammock functions  introduced by S.~Brenner. When dealing
with cluster categories (and cluster-tilted algebras), one should look at
a corresponding class of functions defined on stable
translation quivers, namely the cluster-additive ones.
We conjecture that the cluster-additive functions on a stable translation
quiver of Dynkin type $\ssize \Bbb A_n, \Bbb D_n, \Bbb E_6, \Bbb E_7, \Bbb E_8$ 
are non-negative linear combinations of cluster-hammock
functions (with index set a tilting set). The present paper 
provides a first study of cluster-additive functions and gives
a proof of the conjecture in the case $\ssize \Bbb A_n$.
\par}
	\bigskip\bigskip
A {\it translation quiver} is of the form $\Gamma = (\Gamma_0,\Gamma_1,\tau)$, where
$(\Gamma_0,\Gamma_1)$ is a locally finite quiver say  with $m_{xy}$ arrows $x \to y$, 
and $\tau\:(\Gamma_0\setminus\Gamma_0^p) \to \Gamma_0$
is an injective function defined on the complement of a subset $\Gamma_0^p \subseteq \Gamma_0$, such that 
for any pair of vertices $y,z\in \Gamma_0,$ with $z\notin\Gamma_0^p$ one has $m_{\tau z,y} = m_{y,z}.$
The vertices in $\Gamma_0^p$ are said to be the {\it projective} vertices,
those not in the image of $\tau$ the {\it injective} vertices. If there are neither
projective nor injective vertices, then $\Gamma$ is said to be {\it stable.} 
A typical example of a translation quiver is the Auslander-Reiten quiver
of a finite-dimensional $k$-algebra $A$, where $k$ is an algebraically closed field. Such an Auslander-Reiten
quiver is equipped with an  additive function on the set of vertices with values in the
set of positive integers, its value at a vertex $x$ is the length of the corresponding $A$-module.
Here, a function $f\:\Gamma_0 \to \Bbb Z$ is said to be {\it additive} provided 
$$
 f(z)+f(\tau z) = \sum\nolimits_{y\in \Gamma_0} m_{yz}f(y), \quad \text{for all\ } z\in 
\Gamma_0\setminus \Gamma_0^p.
$$
The importance of dealing with additive functions on translation quivers is well-known 
since a long time, of particular relevance have been the hammock functions introduced
by Brenner [Br], see also [RV]; the hammock functions for the vertices of the 
translations quivers of the form $\Gamma = \Bbb Z\Delta$ with $\Delta$ a Dynkin
diagram $\Bbb A_n, \Bbb D_n, \Bbb E_6, \Bbb E_7, \Bbb E_8$ 
have been displayed already by Gabriel [G] in 1980. 

The present note is concerned with combinatorial features of cluster categories 
(introduced by Buan,  Marsh,  
Reineke,  Reiten,  Todorov [B-T])
and cluster-tilted algebras (introduced by Buan, Marsh, Reiten [BMR]),
and for simplicity we again will assume that we work over an 
algebraically closed field. 
The cluster categories are triangulated categories with Auslander-Reiten triangles, thus we
may consider the corresponding Auslander-Reiten quivers: these are now stable translation 
quivers. Thus, let $\Gamma$ be a stable translation quiver. 
Instead of looking at additive 
functions on $\Gamma$, 
we now will be interested in what we call cluster-additive functions. 

We use the following notation: Any integer $z$ can be written as $z = z^+-z^-$ with non-negative
integers $z^+,z^-$ such that $z^+z^- = 0$ (thus $z^+ = \max\{z,0\}$ and
$z^- = \max\{-z,0\}$). 
A function $f\:\Gamma_0 \to \Bbb Z$ is said to be {\it cluster-additive} on $\Gamma$  provided
$$
 f(z)+f(\tau z) = \sum\nolimits_{y\in \Gamma_0} m_{yz}f(y)^+, \quad \text{for all\ } z\in \Gamma_0.
$$
	\medskip
The sum of two cluster-additive functions usually is not cluster-additive.
Theorem 1 will provide a criterion for such sums to be cluster-additive
again: {\it If $f,g$ are cluster-additive on $\Gamma$, then $f+g$ is 
cluster-additive if and only if $f$ and $g$ are compatible} (this means that
$f(x)g(x) \ge 0$ for all vertices $x$).
Theorem 2 shows that 
{\it the difference $f-g$ of cluster-additive functions 
$f,g$ is cluster-additive  if and only if $g \le f$} (this means that
$g(x)^+ \le f(x)^+$ and $g(x)^- \le f(x)^-$ for all vertices $x$).

	\medskip
The remaining parts of the paper will deal with translation quivers
related to those of the form $\Bbb Z\Delta$ where $\Delta$ is a finite
directed quiver. 
Recall that any locally finite directed quiver $\Delta$  gives rise to a stable 
translation quiver $\Bbb Z\Delta$ with vertex set $\Delta\times \Bbb Z$, with arrows 
$(\alpha,i)\:(\xi,i) \to (\eta,i)$ and $(\alpha^*,i)\:(\eta,i) \to (\xi,i+1)$ 
for any arrow $\alpha\:\xi \to \eta$
in $\Delta$ and with translation $( \xi,i) \mapsto (\xi,i-1)$. 
Theorem 3 asserts that {\it a cluster-additive function on $\Bbb Z\Delta$
with $\Delta$ a finite directed quiver is uniquely determined by its
values on a slice and that these values are arbitrary integers.}
Thus, if $\Delta$ has $n$ vertices, we may identify in this way 
the set of cluster-additive functions on $\Gamma$ with the set $\Bbb Z^n$;
but note that this is just a set-theoretical bijection!

Our main interest lies in the translation quivers $\Bbb Z\Delta$ where
$\Delta$ is a simply laced Dynkin-quiver, thus of type
$\Bbb A_n, \Bbb D_n, \Bbb E_6, \Bbb E_7,$ or $\Bbb E_8$.
Theorem 4 asserts that for $\Delta$ of type $\Bbb A_n$, {\it any
cluster-additive function on $\Bbb Z\Delta$ is a non-negative
linear combination of cluster-hammock functions} (they are
introduced in section 5). We conjecture that the same assertion holds
for all Dynkin cases. This would be an analog of an old theorem 
of Butler [Bu] which asserts that for a representation-finite
algebra $A$, the additive functions on the Auslander-Reiten quiver of $A$
are the linear combinations of the hammock functions. 
	
Cluster-additive categories arise naturally in the context of cluster categories
and cluster-tilted algebras (see section  10), thus one may be tempted to
focus the attention to cluster-additive functions on stable translation
quivers $\Gamma$ such as the Auslander-Reiten quiver of a cluster category,
a typical example is 
$\Bbb Z\Delta/F$ where $\Delta$ is a Dynkin quiver and $F = \tau^{-1}[1].$
It may come as a surprise that instead of looking 
at $\Bbb Z\Delta/F$, we prefer to
consider cluster-additive functions on its cover $\Bbb Z\Delta$. After all,
every cluster-additive function on  $\Bbb Z\Delta/F$ lifts to a cluster-additive
function on $\Bbb Z/\Delta$, thus we deal with a setting which on a fist sight
appears to be more general. But we conjecture that all the cluster-additive 
functions on $\Bbb Z/\Delta$ actually 
are $F$-invariant, so that we would get the shift $F$ for free.

The experienced reader will observe that the cluster-additive functions exhibit
a lot of typical features known in cluster theory (as started by Fomin and
Zelevinsky and developed further by a large number of mathematicians):
that negative numbers arise only seldom, that they have to be ignored in
some calculations, that there is a playing field which concerns only non-negative
numbers, and if the ball leaves the field, it is bounced back immediately ... .
	\medskip
{\bf Acknowledgment.} The paper was written during a stay at the Hausdorff Research
Institute for Mathematics, Bonn, January - April 2011, 
and is inspired by a number of lectures on the 
combinatoric of cluster categories. 
    \bigskip\bigskip
{\bf 1. Preliminaries.}
	\medskip
Let $\Gamma$ be a stable translation quiver. We compare additivity and cluster-additivity and look for the image of a cluster-additive function.
	\medskip
{\bf (1)} {\it A function  $f$ on $\Gamma_0$ with values in $\Bbb N_0$ is cluster-additive
if and only if it is additive.} 
	\medskip
{\bf (2)} {\it If $\Gamma$ is connected and
$\Gamma_1$ is not empty, then any function
 $f\:\Gamma_0 \to \Bbb Z$ which is both additive and
cluster-additive takes values in $\Bbb N_0.$}
	\medskip
Proof: (1) is obvious. For the proof of (2), observe that
a connected stable translation quiver with at least one arrow
has the property that any vertex $y_0$ 
is starting point of an arrow, say $y_0\to z$. Now
$$
 f(z)+f(\tau z) = \sum_y m_{yz}f(y)^+ = \sum_y m_{yz}f(y)
$$ 
implies that $\sum_y m_{yz}(f(y)^+ - f(y)) = 0.$ However we have
$f(y)^+-f(y) \ge 0$ for all $y$. This shows that for  $m_{yz} \neq 0$
we must have $f(y)^+ = f(y).$ Since $m_{y_0,z} \neq 0$, we see
that $f(y_0)^+ = f(y_0).$ Thus $f(y_0) \ge 0.$ 
	\medskip
In the case $\Gamma = \Bbb Z\Bbb A_1$ with vertices $x_i$  $(i\in \Bbb Z)$
such that $\tau x_i = x_{i-1}$, the additive functions are the 
cluster-additive functions, and these are the functions 
of the form $f(x_i) = (-1)^ia$,
where $a$ is a fixed integer. 

	\bigskip
{\bf (3)} {\it Let $f$ be cluster-additive. Let $f(z) < 0$. Then $f(\tau z) \ge -f(z) > 0$.} 
     \medskip
Proof: By definition, $f(\tau z)+ f(z) $ is a sum of positive numbers, thus non-negative. 
	\medskip
This shows:
	\medskip
{\bf (4)} {\it Any cluster-additive function with only non-positive values is the zero function.}  
	\bigskip
{\bf (5)} {\it Let $\Gamma = \Bbb Z\Delta$  with $\Delta$ of Dynkin type.
Any cluster-additive function on $\Gamma$ 
with only non-negative values is the zero function.} 
	\medskip
Proof. Let $f$ be cluster-additive on $\Gamma$ with only non-negative values.
Then $f$  is additive, but according to
[HPR] any additive function on $\Gamma$ with only non-negative values is the zero function. 
	\medskip
It follows from (3) and (5) that there are many stable translation quivers without
non-zero cluster-additive functions. For example, {\it if $\Delta$ is a Dynkin quiver,
then the only cluster-additive function $f$ on $\Gamma = \Bbb Z\Delta/\tau$ is the
zero function.} Namely, (3) asserts that $f$ only takes non-zero values, thus $f$
is additive. Therefore $f$ gives rise to an additive function on $\Bbb Z\Delta$ with
non-negative values. According to (5) this implies that $f$ is the zero function. 

	\bigskip\bigskip
{\bf 2. Sums of cluster-additive functions.}
	\medskip
The sum of two cluster-additive functions usually will not be cluster-additive, a typical example is the following:

	\medskip
{\bf Example.} Let $\Gamma = \Bbb ZA_2.$
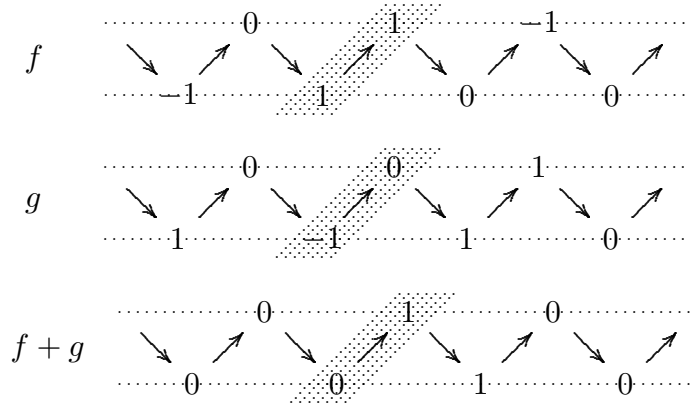
$$
{\beginpicture
\setcoordinatesystem units <.8cm,.8cm>
\put{\beginpicture
\multiput{$0$}  at  -1 1  2 0  4 0 /
\multiput{$-1$}  at -2 0  3 1   /
\multiput{$1$} at 0 0  1 1 /
\setshadegrid span <.5mm>

\hshade -0.3 -.8 0.2  1.3 .8 1.8 /
\arr{-2.7 0.7}{-2.3 0.3}
\arr{-1.7 0.3}{-1.3 0.7}
\arr{-.7 0.7}{-.3 0.3}
\arr{0.3 0.3}{.7 0.7}
\arr{1.3 0.7}{1.7 0.3}
\arr{2.3 0.3}{2.7 0.7}
\arr{3.3 0.7}{3.7 0.3}
\arr{4.3 0.3}{4.7 0.7}
\setdots <1mm>
\plot -3 0  5 0 /
\plot -3 1  5 1 /
\put{$f$} at -4 .5
\endpicture} at 0 0
\put{\beginpicture
\multiput{$0$}  at  -1 1  1 1    4 0 /
\multiput{$-1$}  at  0 0    /
\multiput{$1$} at -2 0  2 0  3 1  /
\setshadegrid span <.5mm>

\hshade -0.3 -.8 0.2  1.3 .8 1.8 /
\arr{-2.7 0.7}{-2.3 0.3}
\arr{-1.7 0.3}{-1.3 0.7}
\arr{-.7 0.7}{-.3 0.3}
\arr{0.3 0.3}{.7 0.7}
\arr{1.3 0.7}{1.7 0.3}
\arr{2.3 0.3}{2.7 0.7}
\arr{3.3 0.7}{3.7 0.3}
\arr{4.3 0.3}{4.7 0.7}
\setdots <1mm>
\plot -3 0  5 0 /
\plot -3 1  5 1 /
\put{$g$} at -4 .5
\endpicture} at 0 -2
\put{\beginpicture
\multiput{$0$}  at  -2 0  -1 1  0 0  3 1  4 0  /
\multiput{$-1$}  at    /
\multiput{$1$} at 1 1 2 0  /
\setshadegrid span <.5mm>

\hshade -0.3 -.8 0.2  1.3 .8 1.8 /
\arr{-2.7 0.7}{-2.3 0.3}
\arr{-1.7 0.3}{-1.3 0.7}
\arr{-.7 0.7}{-.3 0.3}
\arr{0.3 0.3}{.7 0.7}
\arr{1.3 0.7}{1.7 0.3}
\arr{2.3 0.3}{2.7 0.7}
\arr{3.3 0.7}{3.7 0.3}
\arr{4.3 0.3}{4.7 0.7}
\setdots <1mm>
\plot -3 0  5 0 /
\plot -3 1  5 1 /
\put{$f+g$} at -4 .5
\endpicture} at 0 -4

\endpicture}
$$
	\bigskip
Two cluster-additive functions $f,g$ on $\Gamma$ are said to be {\it compatible}
provided $f(x)g(x) \ge 0$ for all vertices $x$. Compatibility can be characterized
in many different ways (the proof is obvious):
	\medskip
{\bf Lemma.} {\it Let $f_1,\dots, f_n$ be cluster-additive functions on $\Gamma$.
The following conditions are equivalent:
\item{\rm (i)} $f_1,\dots,f_n$ are pairwise compatible.
\item{\rm (ii)} If $f_i(x) < 0$ for some index $i$ and some vertex $x$, then $f_j(x) \le 0$
 for $1\le j \le n$.
\item{\rm (iii)} If $f_i(x) > 0$ for some index $i$ and some vertex $x$, then $f_j(x) \ge 0$
 for $1\le j \le n$.
\item{\rm (iv)} Given a pair $i\neq j$, there is no vertex $x$ with $f_i(x) < 0$
 and $f_j(x) > 0.$}
	
	\medskip
{\bf Theorem 1.} {\it Let $f_1,\dots,f_a$ be cluster-additive functions on $\Gamma$.  
Then $\sum f_i$ is cluster-additive if and only if the functions are pairwise compatible.}
	\medskip
Before we start with the proof, let us isolate a decisive property of the operator $z \mapsto z^+$.
	\medskip
{\bf Lemma.} {\it Let $a_1,\dots, a_n$ be integers. Then } 
	\smallskip
(a) {\it  $(\sum_i a_i)^+ \le \sum_i a_i^+.$}
	\smallskip
(b) {\it Equality holds if and only if either all $a_i$ are non-negative or all are non-positive.}
	\medskip
Proof: Let $a_i \ge 0$ for $1\le i \le m$ and $a_i \le 0$ for $m+1 \le i \le n$, let $a = 
\sum_{i=1}^n a_i$. 
Then $\sum_{i=1}^n a_i^+ = \sum_{i=1}^m a_i \ge \sum_{i=1}^n a_i = a$ and therefore
$\sum_{i=1}^n a_i^+ \ge a^+.$ If we have equality, and $a \ge 0,$ then $0 = a - \sum_{i=1}^n a_i^+ 
= \sum_{i=m+1}^n (-a_i)$ shows that these $a_i = 0,$ since all 
$-a_i$ are non-negative for $m+1\le i \le n.$ In this case all the $a_i$ are non-negative.
If $a \le 0$, then $\sum_{i=1}^m a_i = 0$ shows that these $a_i = 0,$ since all $a_i$ are non-negative
for $1\le i \le m.$ In this case, all $a_i$ are non-positive.

Also the converse holds: If all $a_i$ are non-negative, then $\sum_{i=1}^n a_i^+ = 
\sum_{i=1}^n a_i = (\sum_{i=1}^n a_i)^+.$ If all $a_i$ are non-positive, then also 
$\sum_{i=1}^n a_i$ is non-positive, and $\sum_{i=1}^n a_i^+ = 0 = (\sum_{i=1}^n a_i)^+.$
	\bigskip
Proof of Theorem 1: If $\Gamma$ is of tree type $\Bbb A_1$, then the assertion is clear. Thus
we can assume that for any vertex $y$ in $\Gamma$, there is a vertex $z$
with $m_{yz} \neq 0.$ 

First let us assume that $f_1,\dots, f_a$ are pairwise compatible and let $f = \sum_i f_i$.
We claim that for all vertices $y$ of $\Gamma$
$$
 f(y)^+ = \sum\nolimits_i f_i(y)^+ \tag$*$
$$

Let $\Cal T$ be the set of vertices $x\in \Gamma_0$ such that $f_i(x) < 0$ for at least one $i$.
If $y\in \Cal T$, then $f_i(t) \le 0$  for all $1\le i \le a$, since we deal with pairwise compatible 
functions. It follows that $f(y) = \sum_i f_i(y) < 0$ and therefore $f(y)^+ = 0.$ 
But also $f_i(y)^+ = 0$ for all $i$, this
yields $(*)$ in case $y\in \Cal T$. 

Now assume $y\notin \Cal T$. Then $f_i(y) \ge 0$ for all $i$, thus
$f(y) = \sum\nolimits_{i} f_i(y) \ge 0,$
therefore
$$
 f(y)^+ = f(y) = \sum\nolimits_{i} f_i(y) = \sum\nolimits_{i} f_i(y)^+,
$$ 
and we see that $(*)$ is satisfied also in this case. 

Now consider some vertex $z$. 
$$
\align
 f(\tau z)+f(z) &= \sum\nolimits_{i} f_i(\tau z) + \sum\nolimits_{i} f_i(z) 
           = \sum\nolimits_{i} \left(f_i(\tau z)+f_i(z)\right) \cr
           &= \sum\nolimits_{i} \left(\sum\nolimits_y m_{yz}f_i(y)^+\right) \cr
           &= \sum\nolimits_{y} m_{yz} \sum\nolimits_{i} f_i(y)^+ \cr
           &= \sum\nolimits_{y} m_{yz}f(y)^+,
\endalign
$$
where we use that all the functions $f_i$ are cluster-additive as well as the equality $(*)$
for all $y$. This shows that $f$ is cluster-additive.  
	\medskip

Now let us assume that $f = \sum f_i$ is cluster-additive. 
Let $z$ be a vertex of $\Gamma$.
Then, as above, we have
$$
\align
 h(\tau z)+h(z) &= \sum\nolimits_{i} f_i(\tau z) + \sum\nolimits_{i} f_i(z) 
         = \sum\nolimits_{i} \left(f_i(x)+f_i(z)\right) \cr
           &= \sum\nolimits_{i} \left(\sum\nolimits_y m_{yz}f_i(y)^+\right) \cr
           &= \sum\nolimits_{y} m_{yz} \sum\nolimits_{i} f_i(y)^+,
\endalign
$$
thus
$$ 
\align
 0 &= f(\tau z)+f(z) - \sum\nolimits_{i} m_{yz}f(y)^+ \cr
   &= \sum\nolimits_{y} m_{yz} \sum\nolimits_{i} f_i(y)^+ - \sum\nolimits_{y} m_{yz}f(y)^+ \cr
   &= \sum\nolimits_{y} m_{yz} \left(\sum\nolimits_{i} f_i(y)^+ - f(y)^+\right).
\endalign
$$
According to assertion (a) of the Lemma, all the brackets in the last line are non-negative,
thus all the summands $m_{yz} \left(\sum\nolimits_{i} f_i(y)^+ - f(y)^+\right)$ are
non-negative. Since their sum is zero, all these summands are zero.

It follows that for any $y$ we have
$$
  \sum\nolimits_{i} f_i(y)^+ = f(y)^+ 
$$
(since there is $z$ with $m_{yz} \neq 0$). According to the assertion (b) of the Lemma, we 
conclude that all the values $f_i(y)$ for $1\le i \le a$ are non-negative or all are non-positive.
But this means that the functions $f_1,\dots,f_a$ are compatible. 
	\bigskip\bigskip	
   
{\bf 3. Subtraction.}
     \medskip
Let us introduce the following 
{\it partial ordering} on the set of cluster-additive functions on $\Gamma$.
If $f,g$ are cluster-additive functions on $\Gamma$, we write $g \le f$
provided $g(x)^+ \le f(x)^+$ as well as $g(x)^- \le f(x)^-$ for all
vertices $x$ of $\Gamma$. 
	\medskip
{\bf Theorem 2.} {\it Let $f,g$ be cluster-additive functions on
$\Gamma$. Then $f-g$ is cluster-additive if and only if $g\le f.$}
	\medskip
Proof. First, let us assume that $g\le f$. We claim that
$$
 (f-g)(x)^+ = f(x)^+ - g(x)^+
$$
for all vertices $x$.
Namely, if $g(x) > 0,$, then $g(x) = g(x)^+ \le f(x)^+,$ and
therefore $f(x) = f(x)^+$, thus $g(x) \le f(x)$ and therefore
$(f-g)(x) = f(x)-g(x) \ge 0$, thus 
$$
 (f-g)(x)^+ = (f-g)(x) = f(x)-g(x) = f(x)^+-g(x)^+.
$$
Also, if  $g(x) < 0,$ then $g(x)^+ = 0,$ and
$0 < - g(x) = g(x)^- \le f(x)^-$, thus $f(x)^+ = 0.$ Also,
$f(x)^- = - f(x)$ and therefore $g(x) \ge f(x),$ thus
$(f-g)(x) = f(x)-g(x) \le 0.$
It follows that 
$$
 (f-g)(x)^+ = 0 = f(x)^+ - g(x)^+.
$$
Finally, if $g(x) = 0,$ then also $g(x)^+ = 0$ and
$$
 (f-g)(x)^+ = f(x)^+ = f(x)^+ - g(x)^+.
$$

Let $z$ be a vertex of $\Gamma$, then
$$
\align
 (f-g)(\tau z) + (f-g)(z) &=
 f(\tau z) - f(z) + g(\tau z) - g(z) \cr
 &= \sum_y m_{yz}f(y)^+ - \sum_y m_{yz}g(y)^+ \cr
 &= \sum_y m_{yz}(f-g)(y)^+.
\endalign
$$
This shows that $f-g$ is cluster-additive.
	\medskip
Conversely, assume that $f-g$ is cluster-additive. Since
the sum $f = (f-g)+g$ of the cluster-additive functions $f, g$ is
cluster-additive, we know by Theorem 1 that $f-g$ and $g$ are
compatible functions, thus
$(f-g)(x)g(x) \ge 0$ for all vertices $x$, thus
$f(x)g(x) \ge g(x)g(x)$ for all $x$. If $g(x) > 0$, then this
implies that $f(x) \ge g(x) > 0$, thus $f(x)^+ \ge g(x)^+$
and $f(x)^- = g(x)^-$. If $g(x) < 0$, then $f(x) \le g(x) < 0$,
therefore $g(x)^- = -g(x) \le -f(x) = f(x)^-$ and
$g(x)^+ = 0 = f(x)^+.$ Of course, if $g(x) = 0,$ then
$g(x)^+ = 0 \le f(x)^+$ and $g(x)^- = 0 \le f(x)^-.$ This shows
that $g \le f.$
	\bigskip\bigskip	
   
{\bf 4. The restriction of cluster-additive functions to a slice.}
     \medskip
We consider now cluster-additive functions on a translation quiver 
$\Gamma = \Bbb Z\Delta$, where $\Delta$ is a finite directed quiver.
The subset $\Delta_0\times\{0\}$ is called a slice of $\Gamma$ (all
the slices are obtained by considering the subsets 
$\eta(\Delta'_0\times\{0\})$, where 
$\eta\:\Bbb Z\Delta' \to \Bbb Z\Delta$ is an isomorphisms of translation quivers).
	\medskip
{\bf Theorem 3.} {\it Let $\Delta$ be a finite directed quiver.
Any function $f\:\Delta_0\times \{0\} \to \Bbb Z$ can be extended uniquely to a cluster-additive
function on $\Bbb Z\Delta$.}
	  \medskip
This may be reformulated as follows:
     \medskip
{\bf Corollary.} {\it
The restriction furnishes a  bijection between the set of
cluster-additive functions on $\Bbb Z\Delta$ and the functions $f\:\Delta_0\times \{0\} \to \Bbb Z$.}
		 \medskip
Proof of Lemma: Let $x$ be a source in $\Delta$. Then $f$ is defined for $(\xi,0)$ and all its direct successors,
thus we use the defining property of a cluster-additive function in order to define $f(\xi,1)$. 
Inductively we define in this way $f(\eta,j)$ for all vertices $\eta$ of $\Delta$ and all $j > 0$.
The dual procedure yields the values $f(\eta,j)$ for $j <0.$
    \medskip

{\bf Remark 1.}
Note that we need here that  $\Delta$  is finite. For example, if $\Delta$ is the linearly ordered
quiver of type $\Bbb A_\infty^\infty,$ then any function $f\:\Bbb Z\Delta \to \Bbb N_0$
which is constant on the slices $\Delta\times\{i\}$  for all $i\in\Bbb Z$, is additive, thus cluster-additive
and of course not determined by the value taken on one of these slices. 
	\medskip
We also may look at $\Gamma = \Bbb Z\Delta$ with $\Delta$ a locally finite
(but not necessarily finite) directed quiver.
A slice $S$ of $\Gamma$ 
 may be said to be {\it generating} provided  we obtain all vertices from $S$  of $\Gamma$ 
using reflections at sinks and at sources.
If $\Delta$ a finite, 
then any slice is generating, but in general not.
If $\Gamma = \Bbb Z A_\infty^\infty,$ then a slice $S$  
is generating if and only if no arrow in  $S$ belongs to an
infinite path. The corollary can be generalized as follows:
{\it Let $S$ be a generating slice. Then the restriction function 
$f \mapsto f|S$ is bijective.}

    \medskip
{\bf Remark 2.} The extension of a function  $f\:\Delta_0\times \{i\} \to \Bbb Z$ to a cluster-additive
function on $\Gamma$ can be achieved by using what one may call cluster-reflections. Given a locally finite
quiver $\Delta$ and a vertex $x$ of $\Delta$, which is a sink or a source, then
the {\it cluster-reflection} $\sigma_x$ maps any function
$f\:\Delta_0 \to \Bbb Z$ to the function $\sigma_x f$ with $(\sigma_x f)(y) = f(y)$ for $y\neq x$
and $(\sigma_x f)(x) = - f(x) + \sum_y m_{xy} f(y)^+$ and $\sigma_x f$ should be considered
as a function on $(\sigma_x\Delta)_0$, where $\sigma_x\Delta$ is obtained from $\Delta$ by
changing the orientation of all the arrows involving $x$ (thus replacing a source by a sink and vice versa).
Starting with 
a source  $x$ of  $\Delta = \Delta\times\{i\}$, then we may identify $\sigma_x\Delta$ with the slice
obtained by deleting $x$ and adding $\tau^{-1}x$; given a function $f\:\Delta\to \Bbb Z$, and looking for
its cluster-addition extension, then we have to use $\sigma_x f$ on the slice $\sigma_x\Delta$.

Altogether we see that  the restrictions of a cluster-additive function on $\Gamma$ to the 
various slices of $\Gamma$ are obtained from each other by a sequence of cluster-reflections. 
    \bigskip\bigskip

{\bf 5. Cluster-hammock functions.}
     \medskip
Here we introduce some basic cluster-additive functions. As before, we
consider a translation quiver $\Gamma = \Bbb Z\Delta,$ 
where $\Delta$ is a finite directed quiver.
	\medskip
Recall the definition of the {\it left hammock function}
 $h'_p$ for a vertex $p$ of $\Gamma$ 
(and that left hammock functions with finite support are called {\it hammock functions}). First, $h'_p(p) = 1$.
Second, if $z$ is not a successor of $p$, then $h'_p(z) = 0.$ Third, assume that $h'_p(y)$ is defined for all proper predecessors $y$ 
of $z$; if there is an arrow $y \to z$ with $h'_p(y) > 0$, then
$$
 h'_p(z) = - h'_p(\tau z) + \sum\nolimits_{y} m_{yz}h'p(y),
$$
otherwise $h'_p(z) = 0.$ 

It is well-known that all the values $h'_p(z)$ are non-negative; the support of $h'_p$
will be denoted by $H_p$.
If $\Delta$ is a Dynkin quiver (thus of type
$\Bbb A_n, \Bbb D_n, \Bbb E_6, \Bbb E_7,$ or $\Bbb E_8$), then 
$H_p$ is finite and there is a unique vertex $\nu p$ with
$h_p(\nu p) \neq 0$ such that any vertex $y$ with $h_p(y)\neq 0$
is a predecessor of $\nu p$; the map $\nu\:\Gamma_0 \to \Gamma_0$ is
called the {\it Nakayama shift} (see [G], where also typical hammock
functions are displayed; but note that in contrast to the definition
given in this paper, but also in [B] and [RV], Gabriel extends
the function $h'_p|H_p$ to an additive function on all of $\Bbb Z\Delta$).
The shift $\nu\tau^{-1}$ is usually denoted
by $[1]$, the shift $\nu\tau^{-2}$ by $F$. 

We insert here that $\Bbb Z\Delta$ may be interpreted as the
Auslander-Reiten quiver of the derived category $D^b(\mod A)$,
where $A$ is the path algebra of the opposite quiver of $\Delta$, see [H]. Given an 
indecomposable $A$-module $X$, we denote by $[X]$ the corresponding
vertex in $\Bbb Z\Delta$.
In this interpretation,
$[1]$ corresponds to the shift functor of the derived category and
$F$ to the functor $[1]\tau^{-1}$ (also denoted by $F$)
which is used in order to define
the corresponding cluster category, see [B-T].

If $\Delta$ is connected
and not one of these Dynkin diagrams, then the support $H_p$ of $h'_p$ is
infinite, for any vertex $p$ of $\Bbb Z\Delta.$ 
	\medskip

For any vertex $x$ of $\Gamma$, 
we now define a cluster-additive function $h_x$ as follows: 
Let $\Cal S$ be any slice containing $x$, let $h_x(x) = -1$ and $h_x(y) = 0$
for $y\neq x$ in $\Cal S$.

According to Theorem 3, we know that {\it $h_x$ extends 
in a unique way to a cluster-additive function $h_x$ on $\Gamma$ and
this extension is independent of the choice of
$\Cal S$.} We call $h_x$ the {\it cluster-hammock function}
for the vertex $x$.

Proof of the independency: 
There is a slice $\Cal S'$ with $x$ the unique sink of $\Cal S'$ and
a slice $\Cal S''$ with $x$ the unique source of $\Cal S''$, 
and all other slices
containing $x$ are obtained from $\Cal S'$ or also $\Cal S''$
by reflections at sinks or sources different
from $x$. The corresponding cluster-reflections $\sigma_y$ do not change the
value $0$.

Note that the proof shows that $h_x(y) = 0$ for all vertices $y\neq x$
which belong to the convex hull of $\Cal S'$ and $\Cal S''$. 
	
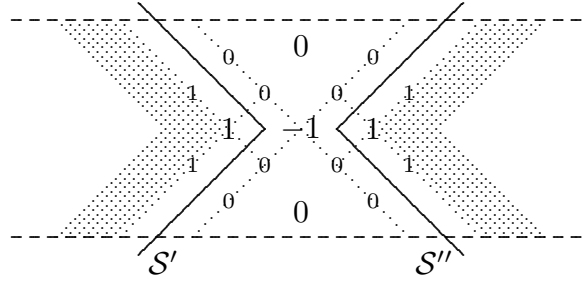
\medskip
$$
{\beginpicture
\setcoordinatesystem units <.4cm,.4cm>
\multiput{} at 0 0  10 6 /
\setdashes <1mm>
\plot -3 0  13 0 /
\plot -3 6  13 6 /

\setdots <1mm>
\plot 0 0  4 4 /
\plot 4 2  0 6 /
\plot 2 0  8 6 /
\plot 2 6  8 0 /
\plot 6 2  7 3  6 4 /
\plot 10 0  7 3  10 6 /

\put{$-1$} at 5 3
\multiput{$\ssize 0$} at 6 4  4 4  6 2  4 2  7 5  7 1  3 1  3 5  /
\multiput{$0$} at  5 .7  5 5.3 /
\multiput{$1$} at  3 3  7 3  /
\multiput{$\ssize 1$} at  2 2  2 4  8 2  8 4   /
\setshadegrid span <.5mm>
\hshade 0 -2 0 <,,z,z> 3 1 3 <,,z,z> 6 -2 0 /
\hshade 0 10 12   <,,z,z> 3 7 9  <,,z,z> 6 10 12 /

\setsolid
\plot .5 -.5  4 3  .5 6.5 /
\plot 9.5 -.5  6 3  9.5 6.5 /
\put{$\Cal S'$} at  1.2 -.7
\put{$\Cal S''$} at 8.7 -.7
\endpicture}
$$

	\medskip
{\bf Lemma.} {\it Let $\Gamma = \Bbb Z\Delta$ with $\Delta$ a Dynkin quiver,
then $h_x$ is $F$-invariant. The support of $h_x$ consists of 
\item{$\bullet$} the $F$-orbit
of $x$ and $h_x$ takes the value $-1$ on these vertices,
as well as 
\item{$\bullet$} the $F$-shifts of the hammock $H_{\tau^{-1}x}$ and
here $h_x$ takes positive values, namely
}
$$
 h_x(y) = h'_{\tau^{-1}x}(y) \quad \text{for} \quad y\in H_{\tau^{-1}x}.
$$
	\medskip
Here is a schematic illustration, where we write $p = \tau^{-1}x$
(the vertical dotted lines mark a fundamental domain for the
action of $F$):
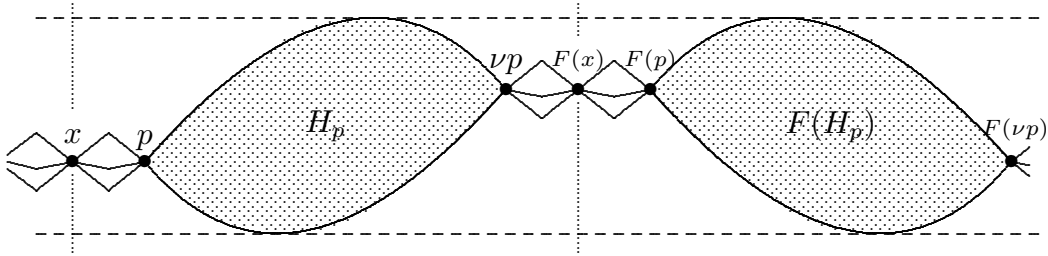
$$
{\beginpicture
\setcoordinatesystem units <.4cm,.4cm>
\multiput{} at 0 0  26 6 /
\setdashes <1mm>
\plot -1 0  27 0 /
\plot -1 6  27 6 /

\setsolid
\multiput{$\bullet$} at  0 2  2 2  12 4  14 4  16 4  26 2 /
\plot 0 2  1 1.2  2 2  1 2.8  0 2  1 1.8  2 2 /
\plot -1.8 1.8  -1 1.2  0 2  -1 2.8  -1.8 2.2 /
\plot -1.8 2  -1 1.8  0 2 /
\plot 12 4  13 3.2  14 4  13 4.8  12 4  13 3.8  14 4 /
\plot 14 4  15 3.2  16 4  15 4.8  14 4  15 3.8  16 4 /
\plot 26 2  26.5 1.6 /
\plot 26 2  26.5 1.9 /
\plot 26 2  26.5 2.4 /
\put{$x$} at 0 2.7 
\put{$p$} at 2 2.6  
\put{$\nu p$} at   12 4.7  
\put{$H_p$} at 7 3  
\put{$\ssize F(x)$} at   14 4.8 
\put{$\ssize F(p)$} at  16 4.8 
\put{$\ssize F(\nu p)$} at    26.1 2.9 
\put{$F(H_p)$} at   21 3 
\setquadratic
\plot 2 2  7.5 5.9  12 4 /
\plot 2 2  6.5 0.1 12 4 /
\plot 2 2  7.5 5.9  12 4 /
\plot 2 2  6.5 0.1 12 4 /
\plot 16 4  20.5 5.9  26 2 /
\plot 16 4  21.5 0.1 26 2 /
\plot 16 4  20.5 5.9  26 2 /
\plot 16 4  21.5 0.1 26 2 /
\setlinear
\setshadegrid span <.5mm>
\vshade 2 2 2 <z,z,,> 3 1 3.1 <z,z,,> 4 0.3 4 <z,z,,> 6 0 5.5 <z,z,,>
 7 0.1 5.9 <z,z,,> 8 0.5 6 <z,z,,>  10 2 5.7 <z,z,,> 11 2.9 5 <z,z,,>
12 4 4 /
\vshade 16 4 4 <z,z,,> 17 2.9 5 <z,z,,> 18 2 5.7 <z,z,,> 20 0.5 6 <z,z,,>
 21 0.1 5.9 <z,z,,> 22 0 5.5 <z,z,,>  24 0.3 4 <z,z,,> 25 1 3.1 <z,z,,>
26 2 2 /
\setdots <.5mm>
\plot 0 -.5   0 1.7 /
\plot 0 3.5   0 6.5 /
\plot 14 -.5   14 3.7 /
\plot 14 5.5   14 6.5 /

\endpicture}
$$

We have mentioned that $h_x$ is $F$-invariant: $h_x(y) = h_x(Fy)$
for all $y\in \Gamma$. But this means also that $h_x = h_{Fx}$.
Thus, when dealing with a set of cluster-hammock functions,
we may restrict to look at those indexed by elements in some fixed
fundamental domain for $F$.
	\medskip
Let us mention a property of the hammock functions $h'_p$ (and of
$h_{\tau p}$) which will be used in the next section. 
If there is a sectional path from $p$ to a vertex $y$, then
$h'_p(y) \ge 1$  (or better: in this case, $h'_p(y)$ is the number
of sectional paths from $p$ to $y$).
	\bigskip
We call a subset $\Cal T$ of $\Bbb Z\Delta$ {\it confined} 
provided there is a slice $\Cal S$ such that $\Cal T$ is contained
in the convex hull of $\Cal S$ and $\tau \Cal S[1]$; note that this is the Auslander-Reiten quiver $\Gamma(A)$ 
of a hereditary algebra $A$ of type $\Delta$, with $\Cal S$ the indecomposable projective $A$-modules, and 
 $\tau \Cal S[1]$ the indecomposable injective $A$-modules. 

We call a subset $\Cal T$ of $\Gamma$ a {\it tilting set} provided we can identify $\Gamma$
as a translation quiver with $D^b(\mod A)$ for some hereditary algebra $A$
such that $\Cal T$ are just the positions of the indecomposable direct
summands of a tilting $A$-module. Subsets of tilting sets are called
{\it partial tilting sets}. 
	\medskip
{\bf Lemma.} {\it  Let $X,Y$ be indecomposable $A$-modules. 
Then $\Ext^1(X,Y) = 0$  if and only if $h_{[Y]}([X]) = 0.$}
     \medskip
Proof: There is the Auslander-Reiten formula $\Ext^1(X,Y) \simeq D\Hom(\tau^-Y,X)$ and
$\dim\Hom(\tau^-Y,X) = h'_{[\tau^{-1}Y]}([X]) = h_[Y]([X]).$

   \medskip
{\bf Corollary.} {\it A subset $\Cal T$ of $\Gamma$ 
is partial tilting if and only if  $\Cal T$ is confined and 
$h_x(x')= 0$  for all pairs $x\neq x'$ in $\Cal T$.}
	
     \bigskip\bigskip
{\bf 6. Non-negative linear combinations of cluster-hammock functions.}
     \medskip
Again, we deal with a translation quiver $\Gamma = \Bbb Z\Delta$,
where $\Delta$ is a finite directed quiver.
	\medskip
{\bf Proposition 1.} {\it Consider a set $h_1,\dots, h_n$ of cluster-hammock functions.
These functions are pairwise compatible if and only if
there is a tilting set $\Cal T$ such that any $h_i$ is of the form $h_x$
with $x\in \Cal T.$ }
	\medskip
Proof: Let $\Cal T$ be a tilting set and $x,x' \in \Cal T$. We have
to show that $h_x,h_{x'}$ are compatible. This is clear if $h_{x} = h_{x'}.$
Thus assume that $h_{x} \neq h_{x'}$, thus $x$ and $x'$ do not belong to
the same $F$-orbit of $\Gamma_0$. Let $h_{x'}(y) < 0$, then $y$ belongs to
the $F$-orbit of $x'$, thus $h_{x}(y) = h_{x}(x') = 0.$ This shows
that $h_{x}(y)h_{x'}(y) = 0$. Similarly, we see: if $h_{x}(y) < 0,$ then
$h_{x}(y)h_{x'}(y) = 0$. For the remaining vertices $y$ we have 
both $h_{x}(y) \ge 0$ and $h_{x'}(y) \ge 0$, thus also $h_{x}(y)h_{x'}(y) \ge 0.$

Conversely, assume that the functions $h_1,\dots,h_n$ are pairwise compatible. 
First, we show that for $h_i\neq h_j$, and $h_j = h_y$ for some vertex $y$,
then $h_i(y) = 0.$ Namely, $h_i(y)h_y(y) = h_i(y)h_j(y) \ge 0$, and $h_y(y) = -1$
shows that $h_i(y) \le 0$. But $h_i(y) < 0$ would imply that $h_i = h_y$,
a contradiction. Thus $h_i(y) = 0.$ 

Now, let $h_1 = h_x$ for some
$x\in \Gamma_0.$ Let $\Cal S$ be the slice in $\Gamma$ such that $\tau^{-2}F^{-1}x$
is the unique source. Let $\Cal S' = \Cal S[1]$, this is the slice
with unique source $\tau^{-1}x.$ Clearly, the convex hull $\Cal F$ of $\Cal S$ and $\Cal S'$
is a fundamental domain for $F$, thus $h_i = h_{x_i}$ for some $x_i\in \Cal F$.
Since $\tau^{-1}x$ is the unique source of $\Cal S'$, we see that $h_x(z) > 0$
for all $z\in \Cal S'$. Assume that some $x_j$ belongs to $\Cal S'$, then
$x \neq x_j$, thus $h_x \neq h_{x_j}$ (since $x,x_j$ belong to the fundamental
domain $\Cal F$ of $F$), but then we know that $h_x(x_j) = 0$, a contradiction.
In this way, we see that all the vertices $x_i$ belong to the convex hull
of $\Cal S$ and $\tau \Cal S' = \tau\Cal S[1]$, thus the set 
$\Cal T = \{x_1,\dots,x_n\}$ is confined. Since also $h_x(x') = 0$
for $x\neq x'$ in $\Cal T,$ we see that $\Cal T$ is a tilting set.

	\bigskip
{\bf Corollary.} 
{\it A linear combination $h = \sum_{x\in \Cal T} n_xh_x$ with positive integers $n_x$ is
cluster-additive if and only if $\Cal T$ is a partial cluster-tilting set.}
	\medskip

Proof. This is a direct consequence of Theorem 1 and Proposition 1.

	\medskip
{\bf Proposition 2.} {\it 
Let $f = \sum_{x\in \Cal T} n_xh_x$ for some tilting set $\Cal T$  and $n_x\in \Bbb N_0$, then
$f(x) = -n_x$ for $x\in \Cal T$ and $f(y) \ge 0$ provided the intersection of $\Cal T$ and the $F$-orbit of $y$ is empty.
Thus}
$$
 f = \sum_{x\in \Cal T} n_xh_x =- \sum_{x\in \Cal T} f(x)h_x =  
\sum_{x\in \Cal T} f(x)^-h_x =  \sum_{x\in \Gamma^0} f(x)^-h_x.
$$
where $\Gamma^0$ is the convex hull of some slice 
$\Cal S$ and $\tau\Cal S[1]$.

	\bigskip
{\bf Conjecture.} {\it Let  $\Gamma = \Bbb Z\Delta$ where $\Delta$ is one of the Dynkin
diagrams $\Bbb A_n, \Bbb D_n, \Bbb E_6, \Bbb E_7, \Bbb E_8$ and let $f$ be cluster-additive on 
$\Gamma$. Then $f$ is a non-negative linear combination of cluster-hammock
functions} (and therefore 
of the form
$$
 \sum\nolimits_{x\in \Cal T} n_xh_x
$$
for a tilting set $\Cal T$ and integers $n_x\in \Bbb N_0$, for all $x\in \Cal T$).
    \medskip

If this conjecture is true, then any cluster-additive function satisfies the following properties:
   \medskip
(a) {\it $f$ is $F$-invariant.}
    \medskip
(b) {\it $\{x\in\Gamma_0\mid f(x) < 0\}$ 
 is the union of the $F$-orbits of a partial tilting set.}
    \medskip
(c) {\it There is a partial tilting set $\Cal T$ with }
$$
 f = \sum\nolimits_{x\in \Cal T} f(x)^-h_x
$$

A proof of the conjecture in the case $\Bbb A_n$ will be given in section 9.
We also note that it is not difficult to exhibit explicitely all the cluster-additive
functions on $\Bbb Z\Delta$, where $\Delta$ is a quiver of type $\Bbb D_4$, thus
verifying the conjecture also in this case.
    
		         \bigskip\bigskip
{\bf 7. The rectangle rule.}
     \bigskip
{\bf Lemma.} {\it Let $f$ be cluster-additive on the following translation quiver
with $s\ge 1,\ t\ge 1$:
$$
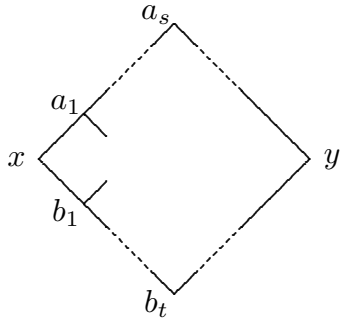
{\beginpicture
\setcoordinatesystem units <.5cm,.5cm>
\plot 1.5 1.5  0 0  1.5 -1.5 /
\plot 2.5 2.5   3 3  3.5 2.5 /
\plot 4.5 1.5 6 0  4.5 -1.5 /
\plot 2.5 -2.5  3 -3  3.5 -2.5 /

\plot 1 1  1.5 0.5 /
\plot 1 -1 1.5 -.5 /
\setdashes <.5mm> 
\plot 1.5 1.5  2.5 2.5 /
\plot 1.5 -1.5  2.5 -2.5 /
\plot 4.5 1.5  3.5 2.5 /
\plot 4.5 -1.5  3.5 -2.5 /
\put{$x$} at -.5 0
\put{$y$} at 6.5 0
\put{$a_1$} at 0.6 1.2 
\put{$a_s$} at 2.6 3.2 
\put{$b_1$} at 0.6 -1.2 
\put{$b_t$} at 2.6 -3.2 
\endpicture}
$$
Then for $y = y(s,t)$ with $f(x) \le 0$, 
we have
$$
 f(y) = f(x)^- +\sum_{1\le i \le s-1} f(a_i)^- + f(a_s)^+ + \sum_{1\le j \le t-1} f(b_j)^- + f(b_t)^+.
$$
In particular, $f(y) \ge f(x)^- \ge 0.$} 
   \medskip
Proof, by induction on $s$ and $t$.

If $s=t = 1$, then $f(y) = f(a_1)^+ + f(b_1)^+ - f(x).$

Now assume that we know the formula for some $s,t$. Let us increase $s$ by $1$, thus we deal with
$$
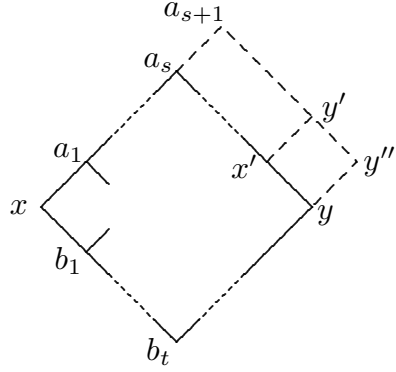
{\beginpicture
\setcoordinatesystem units <.5cm,.5cm>
\plot 1.5 1.5  0 0  1.5 -1.5 /
\plot 2.5 2.5   3 3  3.5 2.5 /
\plot 4.5 1.5 6 0  4.5 -1.5 /
\plot 2.5 -2.5  3 -3  3.5 -2.5 /

\plot 1 1  1.5 0.5 /
\plot 1 -1 1.5 -.5 /
\setdashes <.5mm> 
\plot 1.5 1.5  2.5 2.5 /
\plot 1.5 -1.5  2.5 -2.5 /
\plot 4.5 1.5  3.5 2.5 /
\plot 4.5 -1.5  3.5 -2.5 /
\put{$x$} at -.5 0
\put{$y$} at 6.3 -.2
\put{$y'$} at 6.5 2.2
\put{$y''$} at 7.5 1
\put{$a_1$} at 0.6 1.2 
\put{$a_s$} at 2.6 3.2 
\put{$a_{s+1}$} at 3.4 4.3 
\put{$b_1$} at 0.6 -1.2 
\put{$b_t$} at 2.6 -3.2 
\setdashes <1mm>
\plot 3 3  4 4   7 1  6 0 /
\plot 5 1  6 2 /
\put{$x'$} at 4.5 .9 
\endpicture}
$$
For $t=1$, we have $x' = a_s$ and $y' = a_{s+1}$,
otherwise $x' = y(s,t-1)$ and $y' = y(s+1,t-1).$

Now, consider first the case $t=1$. Then (since $f(y) \ge 0$):
$$
\align
 f(y'') &= f(y')^+ + f(y) - f(x') \cr
        &= 
          f(a_{s+1})^+ + f(x)^- +\sum_{1\le i \le s-1} f(a_i)^- + f(a_s)^+ 
                            + \sum_{1\le j \le t-1} f(b_j)^- + f(b_t)^+
        - f(a_s) \cr
        & = f(a_{s+1})^+ + f(x)^- +\sum_{1\le i \le s} f(a_i)^- 
                            + \sum_{1\le j \le t-1} f(b_j)^- + f(b_t)^+
\endalign
$$
where the last equality comes from $f(a_s)^+ - f(a_s) = f(a_s)^-$.

Second, let $t \ge 2$. Then both $f(y)\ge 0,\ f(y') \ge 0$, thus
$$
\align
 f(y'') &= f(y') + f(y) - f(x') \cr
        &=\quad  f(x)^- \ +\ \sum_{1\le i \le s} f(a_i)^- \ + f(a_{s+1})^+ 
                            + \sum_{1\le j \le t-2} f(b_j)^- + f(b_{t-1})^+ \cr
        &\quad\; +  f(x)^- +\sum_{1\le i \le s-1} f(a_i)^- + f(a_s)^+ 
                          \ \ + \sum_{1\le j \le t-1} f(b_j)^- + f(b_t)^+ \cr
        &\quad\; -  f(x)^- -\sum_{1\le i \le s-1} f(a_i)^- - f(a_s)^+ 
                        \  \  - \sum_{1\le j \le t-2} f(b_j)^- - f(b_{t-1})^+ \cr
        &=\quad  f(x)^- \ +\ \sum_{1\le i \le s} f(a_i)^- \ + f(a_{s+1})^+  
                        + \sum_{1\le j \le t-1} f(b_j)^- + f(b_t)^+,
\endalign
$$
as we want. 

By symmetry, the same argument works, if we increase $t$ instead of $s$. This completes the
proof.
	\bigskip
{\bf Extended version.} {\it Let $f$ be cluster-additive on the following translation quiver
with $s\ge 1,\ t\ge 1$:
$$
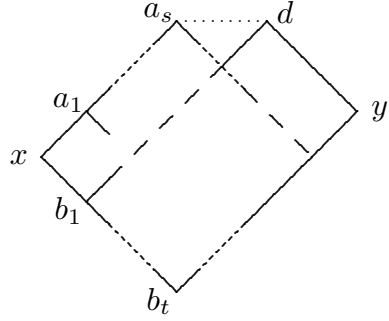
{\beginpicture
\setcoordinatesystem units <.5cm,.5cm>
\plot 1.5 1.5  0 0  1.5 -1.5 /
\plot 2.5 2.5   3 3  3.5 2.5 /
\plot 6 0  4.5 -1.5 /
\plot 2.5 -2.5  3 -3  3.5 -2.5 /

\plot 1 1  1.5 0.5 /
\plot 1 -1 1.5 -.5 /

\setdashes <.5mm>

\plot 1.5 1.5  2.5 2.5 /
\plot 1.5 -1.5  2.5 -2.5 /
\plot 4.5 1.5  3.5 2.5 /
\plot 4.5 -1.5  3.5 -2.5 /
\put{$x$} at -.5 0
\put{$y$} at 7.5 1
\put{$a_1$} at 0.6 1.2 
\put{$a_s$} at 2.6 3.2 
\put{$b_1$} at 0.6 -1.2 
\put{$b_t$} at 2.6 -3.2 

\setsolid
\plot 4 2  5 3  7 1  6 0 /
\setdots <1mm>
\plot 3 3  5 3 /
\put{$d$} at 5.4 3.2 
\setdashes <2mm>
\plot 1.5 -.5  4 2 /

\plot 4.5 1.5  6 0 /

\endpicture}
$$
Then for $y = y(s+1,t)$ with $f(x) \le 0$, 
we have}
$$
 f(y) = f(x)^- +\sum_{1\le i \le s-1} f(a_i)^- + \sum_{1\le j \le t-1} f(b_j)^- + f(b_t)^+.
$$
	\medskip
Proof: We add a vertex $a_{s+1}$ and arrows $a_s \to a_{s+1}$ and $a_{s+1} \to d$,
so that we obtain a rectangle. Also, we extend $f$ to be defined on the rectangle
by setting $f(a_s) = 0.$ Then the extended function satisfies the
cluster-additivity condition on all the meshes of the rectangle and we can apply the
lemma.
$$
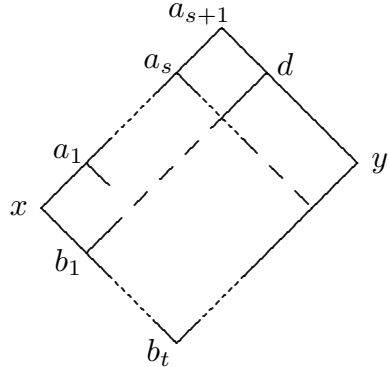
{\beginpicture
\setcoordinatesystem units <.5cm,.5cm>
\plot 1.5 1.5  0 0  1.5 -1.5 /
\plot 2.5 2.5   3 3  3.5 2.5 /
\plot  6 0  4.5 -1.5 /
\plot 2.5 -2.5  3 -3  3.5 -2.5 /

\plot 1 1  1.5 0.5 /
\plot 1 -1 1.5 -.5 /
\setdashes <.5mm> 
\plot 1.5 1.5  2.5 2.5 /
\plot 1.5 -1.5  2.5 -2.5 /
\plot 4.5 1.5  3.5 2.5 /
\plot 4.5 -1.5  3.5 -2.5 /
\put{$x$} at -.5 0
\put{$y$} at 7.5 1
\put{$a_1$} at 0.6 1.2 
\put{$a_s$} at 2.6 3.2 
\put{$b_1$} at 0.6 -1.2 
\put{$b_t$} at 2.6 -3.2 

\setsolid
\plot 4 2  5 3  7 1  6 0 /
\put{$d$} at 5.4 3.2 
\setdashes <2mm>
\plot 1.5 -.5  4 2 /

\plot 4.5 1.5  6 0 /

\setsolid 
\plot 3 3  4 4  5 3 /

\put{$a_{s+1}$} at 3.5 4.3
\endpicture}
$$

	\bigskip
There is also a corresponding double extended version for dealing with $\Bbb Z\Delta$
where $\Delta$ is of type $\Bbb A_{s+t+1}.$
	\medskip
{\bf Double extended version.} {\it Let $f$ be cluster-additive on the following translation quiver
with $s\ge 1,\ t\ge 1$:
$$
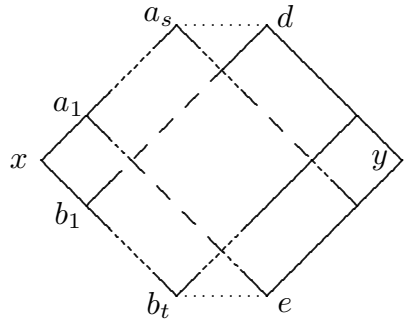
{\beginpicture
\setcoordinatesystem units <.5cm,.5cm>
\plot 1.5 1.5  0 0  1.5 -1.5 /
\plot 2.5 2.5   3 3  3.5 2.5 /
\plot 6 0  4.5 -1.5 /
\plot 2.5 -2.5  3 -3  3.5 -2.5 /

\plot 1 1  1.5 0.5 /
\plot 1 -1 1.5 -.5 /
\setdashes <.5mm> 
\plot 1.5 1.5  2.5 2.5 /
\plot 1.5 -1.5  2.5 -2.5 /
\plot 4.5 1.5  3.5 2.5 /
\plot 4.5 -1.5  3.5 -2.5 /
\plot 1.5 0.5  2 0 /

\plot 6 0  6.5 -0.5  /
\plot 4 -2  4.5 -2.5 /

\put{$x$} at -.5 0
\put{$y$} at 7.5 0
\put{$a_1$} at 0.6 1.2 
\put{$a_s$} at 2.6 3.2 
\put{$b_1$} at 0.6 -1.2 
\put{$b_t$} at 2.6 -3.2 

\setsolid
\plot 4 2  5 3  7 1  6 0 /
\plot 4.5 -2.5 5 -3  8 0  7 1 /
\plot 6.5 -0.5  7 -1 /

\setdots <1mm>
\plot 3 3  5 3 /
\plot 3 -3  5 -3 /
\put{$d$} at 5.4 3.2 
\put{$e$} at 5.4 -3.2 
\setdashes <2mm>
\plot 1.5 -.5  4 2 /
\plot 2 0  4 -2 /
\plot 4.5 1.5  6 0 /

\endpicture}
$$
Then for $y = y(s+1,t+1)$ with $f(x) \le 0$, 
we have}
$$
 f(y) = f(x)^- +\sum_{1\le i \le s-1} f(a_i)^- + \sum_{1\le j \le t-1} f(b_j)^-.
$$

		 \bigskip\bigskip
{\bf 8. Wings.}
	\medskip
Let $s\ge 0, t\ge 1$, let $y$ be a wing vertex of rank $s+t+1$, say with
sectional paths
$$
 p[1] \to p[2] \to \cdots \to p[s+t+1] = y,\qquad 
y =  [s+t+1]q \to \cdots \to  [2]q \to q[1].
$$

{\bf Lemma.} {\it Assume that 
$$
 f(p[s]) \le 0,\quad f(p[s+i]) \ge 0, \text{ for }  1\le i \le t,\quad f(p[s+t+1]) \le 0.
$$
Then 
$$
 f([t]q) = \min_{1 \le i \le t} f(p[s+i]).
$$ 
Also, $f$ is non-negative on all vertices between $p[s+1]$ and
$[1+t]q$ different from $y$.} 
	\medskip
Here is a sketch which exhibits the vertices in question in case $s \ge 1$:
$$
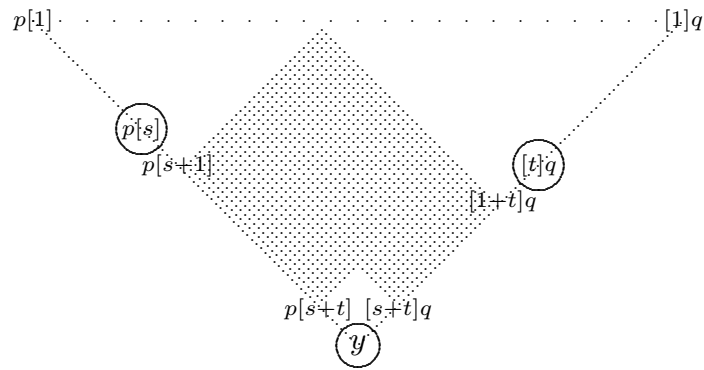
{\beginpicture
\setcoordinatesystem units <.4cm,.4cm>
\multiput{} at 0 0  18 -9 /
\put{$\ssize p[1]$} at 0 0 
\put{$\ssize p[s]$} at 3 -3
\put{$\ssize p[s+1]$} at 4 -4 
\put{$\ssize p[s+t]$} at 7.9 -8
\put{$y$} at 9 -9
\put{$\ssize [s+t]q$} at 10.1 -8
\put{$\ssize [1+t]q$} at 13 -5
\put{$\ssize [1]q$} at 18 0

\put{$\ssize [t]q$} at 14 -4 

\setdots <.7mm>
\plot 0 0  9 -9  18 0 /
\setdots <3mm> 
\plot 0 0  18 0 /
\setsolid
\circulararc 360 degrees from 9 -9.6 center at 9 -9
\circulararc 360 degrees from 3 -3.7 center at 3 -3
\circulararc 360 degrees from 14 -4.7 center at 14 -4

\setdots <2mm> 

\setshadegrid span <.5mm>

\vshade 4 -4 -4 <z,z,,>  8 -8 0  <z,z,,> 9 -7 -1  <z,z,,> 10 -8 -2   <z,z,,> 13 -5 -5 /

\endpicture}
$$
The case $s = 0$ looks as follows:
$$
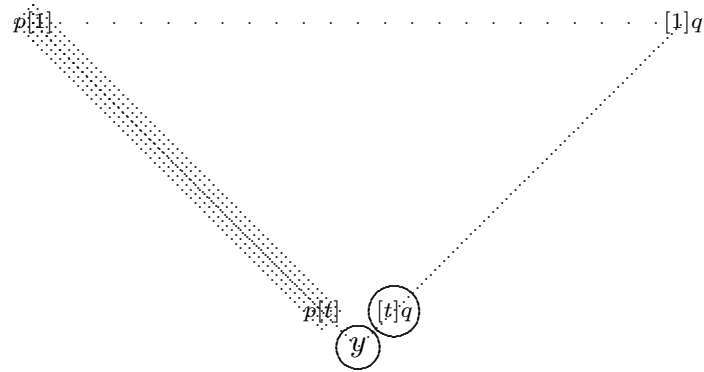
{\beginpicture
\setcoordinatesystem units <.4cm,.4cm>
\multiput{} at 0 0  18 -9 /
\put{$\ssize p[1]$} at 0 0 
\put{$y$} at 9 -9
\put{$\ssize p[t]$} at 8 -8
\put{$\ssize [t]q$} at 10 -8
\put{$\ssize [1]q$} at 18 0

\setdots <.7mm>
\plot 0 0  9 -9  18 0 /
\setdots <3mm> 
\plot 0 0  18 0 /
\setsolid
\circulararc 360 degrees from 9 -9.6 center at 9 -9
\circulararc 360 degrees from 10 -8.7 center at 10 -8

\setshadegrid span <.5mm>

\vshade -0.7 0 0 <z,z,,>  0 -0.7 0.7  <z,z,,> 8 -8.7 -7.3  
<z,z,,> 8.7 -8 -8 /

\endpicture}
$$
	\medskip
Proof: Let us use the following labels for the relevant vertices
of the wing
$$
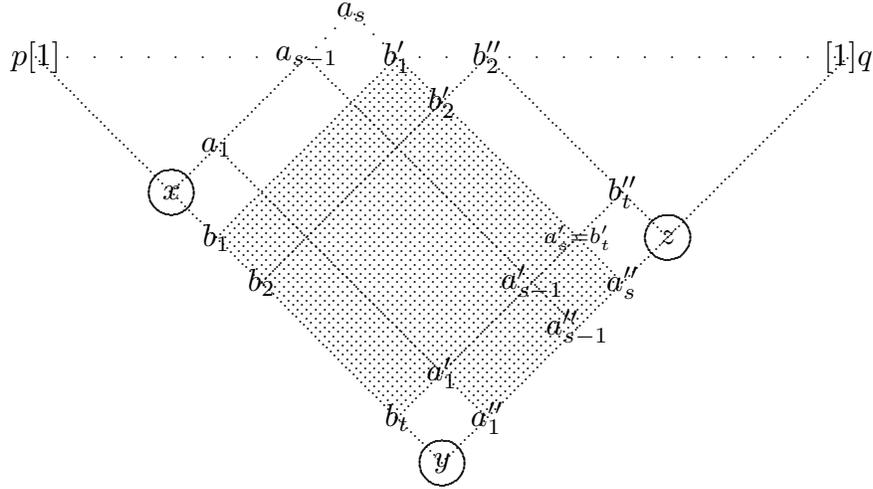
{\beginpicture
\setcoordinatesystem units <.5cm,.5cm>
\multiput{} at 0 0  18 -9 /
\put{$p[1]$} at 0 0 
\put{$x$} at 3 -3
\put{$b_1$} at 4 -4 
\put{$b_2$} at 5 -5 
\put{$b_t$} at 8 -8
\put{$y$} at 9 -9

\put{$a_1''$} at 10 -8
\put{$a_{s-1}''$} at 12 -6  
\put{$a_s''$} at 13 -5

\put{$a_1'$} at 9 -7
\put{$a_{s-1}'$} at 11 -5
\put{$\ssize a_s'=b_t'$} at 12 -4

\put{$[1]q$} at 18 0
\put{$a_1$} at 4 -2
\put{$a_{s-1}$} at 6 0
\put{$a_s$} at 7 1
\put{$b_1'$} at 8 0
\put{$b_2'$} at 9 -1 

\put{$b_2''$} at 10 0
\put{$b_t''$} at 13 -3

\put{$z$} at 14 -4 

\setdots <.6mm>
\plot 0 0  9 -9  18 0 /
\plot 3 -3  6 0  12 -6 /
\plot 4 -4  8 0  13 -5 /
\plot 8 -8  13 -3  14 -4 /
\plot 4 -2  10 -8 /
\plot 5 -5  10 0  13 -3 /

\setdots <3mm> 
\plot 0.7 0   5.5 0 /
\plot 8.5 0  9.5 0 /
\plot 10.5 0  18 0 /

\setsolid
\circulararc 360 degrees from 9 -9.5 center at 9 -9
\circulararc 360 degrees from 3 -3.5 center at 3 -3
\circulararc 360 degrees from 14 -4.5 center at 14 -4

\setdots <1.5mm> 
\plot 6 0  7 1  8 0   /

\setshadegrid span <.5mm>

\vshade 4 -4 -4 <z,z,,>  8 -8 0  <z,z,,> 9 -7 -1  <z,z,,> 10 -8 -2   <z,z,,> 13 -5 -5 /

\endpicture}
$$
In particular
$$
\gather
 x = p[s],\ z = [t]q,\ b_i = p[s+i],\ a_i = \tau^{-i}p[s-i].
\endgather
$$
Note that we have added the vertex $a_s$ 
to the wing, (with additional arrows $a_{s-1} \to a_s$ and
$a_s \to b_1'$)
 and we put $f(a_s) = 0,$ as in the proof of the extended rectangle rule.
	\medskip
Using the new labels, the assumptions read:
$$
 f(x) \le 0,\quad f(b_i) \ge 0,\ \text{ for }\  1\le i \le t,\quad f(y) \le 0
$$
and the assertion is that $f$ is non-negative on the shaded area (the 
vertices between $b_1$ and $a_s''$ different from $y$) 
and that 
$$
 f(z) = - \min(f(b_i)\mid 1\le i \le t).
$$

The rectangle rule asserts that $f$ is bounded below by $f(x)^-$
on the rectangle between $\tau^{-1}x$ and $a_s' = b_t'$. By assumption, $f$ is
non-negative on the vertices $b_1,\dots, b_t$. Thus,
concerning the non-negativity assertion, it remains to show that
$f$ is non-negative on the vertices $a_1'',\dots, a_s''$.

The rectangle rule asserts that
$$
 f(a_i') = f(x)^- + \sum_{j=1}^{i-1}f(a_j)^-+f(a_i)^+ +\sum_{j=1}^{t-1}f(b_j)^- +f(b_t)^+.
$$
Since $f(y) \le 0$ and $f(a_1') \ge 0$, we have $f(a_1'') 
= f(a_1') - f(b_t) \ge 0$.
Assume by induction that we know that 
$f(a_i'') = f(a_i') - f(b_t) \ge 0$, then
we get
$$
\align
 f(a_{i+1}'') + f(a_i') &= f(a_i'')^+ + f(a_{i+1}')^+ \cr
   &= f(a_i'') + f(a_{i+1}') \cr
   &= f(a_i') - f(b_t) + f(a_{i+1}') 
\endalign
$$
and therefore 
$$
 f(a_{i+1}'') = f(a_{i+1}') -  f(b_t).
$$ 
By the rectangle rule for $a_{i+1}'$ we see that $f(a_{i+1}') -  f(b_t)\ge 0$
provided $i+1\le s.$

	\medskip
It remains to calculate the value $f(z).$

Using induction on $i$, we show that 
$$
 f(b_i'') = f(b_i) - \min(f(b_j)\mid 1\le j < i)
$$
for $i\ge 2.$

The rectangle rule for $b_i'$ yields
$$
\align
 f(b_i') & = f(x)^- + \sum_{j=1}^{s-1}f(a_j)^-+f(a_s)^+ +\sum_{j=1}^{i-1}f(b_j)^- +f(b_i)^+ \cr
  & = f(x)^- + \sum_{j=1}^{s-1}f(a_j)^- +f(b_i),
\endalign
$$
since $f(a_s) = 0$ and all $f(b_j) \ge 0$.
Similarly, for $b_{i+1}'$ we get:
$$
 f(b_{i+1}')  = f(x)^- + \sum_{j=1}^{s-1}f(a_j)^- +f(b_{i+1}),
$$
thus
$$
  f(b_{i+1}') -  f(b_{i}') = f(b_{i+1}) - f(b_{i}).
$$
For $i = 2$, we have
$$
 f(b_2'') = f(b_2') - f(b_1'),
$$
since $f(b_2') \ge 0$, thus
$$
 f(b_2'') = f(b_2') - f(b_1') = f(b_2) - f(b_1) = f(b_2) - 
  \min(f(b_j)\mid 1\le j < 2),
$$
as we have claimed. 

Similarly, we have for all $i\ge 2$ 
$$
\align
  f(b_{i+1}'') &=  f(b_{i}'')^+ + f(b_{i+1}') -  f(b_{i}') \cr
   &= f(b_{i}'')^+ + f(b_{i+1}) - f(b_{i}).
\endalign
$$
By induction, we may assume that 
$$
 f(b_i'') = f(b_i) - \min(f(b_j)\mid 1\le j < i),
$$
and we have to distinguish two cases:

First, assume that 
$f(b_i'') \le 0.$ Then $f(b_i'')^+ = 0$ and 
$f(b_i) \le \min(f(b_j)\mid 1\le j < i),$ so that
$\min(f(b_j)\mid 1\le j \le i) = f(b_i)$. Thus 
$$
\align
  f(b_{i+1}'') &= f(b_{i}'')^+ + f(b_{i+1}) - f(b_{i}) \cr
               &= 0 + f(b_{i+1}) - \min(f(b_j)\mid 1\le j \le i),
\endalign
$$
as we want to show. 

In the second case, 
$f(b_i'') \ge 0,$ thus $f(b_i'')^+ = f(b_i'')$ and
$f(b_i) \ge \min(f(b_j)\mid 1\le j < i),$ so that
$\min(f(b_j)\mid 1\le j \le i) = \min(f(b_j)\mid 1\le j < i).$
Thus
$$
\align
  f(b_{i+1}'') &= f(b_{i}'')^+ + f(b_{i+1}) - f(b_{i}) \cr
               &= f(b_i) - \min(f(b_j)\mid 1\le j < i) 
 + f(b_{i+1}) - f(b_{i}) \cr
               &= \min(f(b_j)\mid 1\le j \le i) 
 + f(b_{i+1}).
\endalign
$$

Thus we see that 
$$
 f(b_t'')^+ = f(b_t)  - \min(f(b_i)\mid 1\le i \le t).
$$
On the other hand, the calculations in the first part of the
proof had shown that $f(a_s'') \ge 0$ and  that
$$
 f(a_s'') - f(a_s') = -f(b_t).
$$
It follows that 
$$
\align
 f(z) &= f(b_t'')^+ + f(a_s'') - f(a_s') =
f(b_t'')^+ + f(a_s'')^+ - f(a_s') 
\cr
 &= f(b_t'')^+ - f(b_t) = f(b_t)  - \min(f(b_i)\mid 1\le i \le t) - f(b_t)\cr
 &=   - \min(f(b_i)\mid 1\le i \le t)
\endalign
$$

This completes the proof.
	\bigskip\bigskip

{\bf 9. The case $\Gamma = \Bbb Z\Bbb A_n$}		 			\medskip
Consider now the case $\Gamma = \Bbb Z\Delta$ with $\Delta$ of
type $\Bbb A_n$. 
	\medskip
{\bf Theorem 4.} {\it Let  $\Gamma = \Bbb Z\Delta$ with $\Delta$ of
type $\Bbb A_n$. Then any cluster-additive function on $\Gamma$ is
a non-negative linear combination of cluster-hammock functions.}
	\bigskip
If $n = 1$, then any cluster-additive function on $\Gamma$ is
a non-negative multiple of one of the two cluster-hammock functions.
Thus, we can assume that $n \ge 2$.
	\medskip
Let $f$ be a cluster-additive function on $\Gamma$.
	\medskip
(1) {\it If $z$ is a vertex of $\Gamma$ with $f(z) \le 0$, then there is
a vertex $z'\neq z$ with $f(z') \le 0$ and a sectional path from $z$ to 
$z'$ or from $z'$ to $z$.}
	\medskip
Proof: Since $n \ge 2$, there is an arrows $a_1 \to a_0 = z$.
Choose $m$ maximal such that there exists a sectional path
$$
 a_m \to \cdots \to a_1 \to a_0 = z.
$$
If $f(a_i) \le 0$ for some $1 \le i \le m$, then let $z' = a_i.$
Otherwise we consider the wing with corners
$$
 p[1] = a_m,\quad z, \quad q[1] = \tau^{-m}a_m.
$$
The wing lemma (with $s=0$) asserts that $f(z) \le 0$ (even $f(z) < 0)$
for $z' = \tau^{-1}a_1$.
	\medskip
(2) {\it If $f(z) < 0$ for some vertex $z$, then $f(z) = f(Fz)$
and }
$$
 f(z)^-h_z \le f.
$$
	\medskip
Proof: According to (1), there is a vertex $y = z'$ with $f(y) \le 0$
and a sectional path from $z$ to 
$y$ or from $y$ to $z$. Up to duality, we can assume that there is
a sectional path from $y$ to $z$ (otherwise we consider instead of
$\Gamma$ the opposite
translation quiver). Also, we can assume that we choose $y$ such that
the path from $y$ to $x$ is of smallest possible length (thus $f$ is
positive on all the vertices between $y$ and $z$). Consider the wing
with corners 
$$
  p[1],\quad y,\quad [1]q,
$$
thus there are sectional paths
$$
 p[1] \to p[2] \to \cdots \to p[m] = y,\qquad 
 y =  [m]q \to \cdots \to  [2]q \to q[1]
$$
and $z$ is one of the vertices $[i]q$ with $1\le j < m$.
Let $s\ge 0$ be maximal with $f(p[s]) \le 0$ and $t=m-s-1$.
We claim that $t\ge 1$ and that $z = [t]q.$

First of all, for $t = 0$, the rectangle rule would imply that
$f([j]q) \ge 0$ for $1\le j < m$, but $z$
is of the form $[j]q$ and $f(z) < 0.$ 

This means that we 
can use the wing lemma, it asserts that 
$f([t]q) = \min (f(p[s+i]\mid 1 \le i \le t$ and that
$f([j]q) \ge 0$ for $1+t \le j \le s+t.$ Since $f(z) < 0$,
with $z$ of the form $[j]q$ and $j\le s+t$, it follows that 
$j\le t$. On the other hand, we know that $f$ is positive
on all vertices between $y$ and $z$, thus we see that $j = t.$
	\medskip

Let $x = p[s],$ $b_i = p[s+i]$, and $z = [t]q$ and note that we have
$f(x) \le 0,\ f(y) \le 0,$
and $f(b_i) \ge 1$ for $1\le i \le t.$ 
This yields the upper wing in the following picture, namely the wing with
corners
$$
 p[1],\ y,\ [1]q.
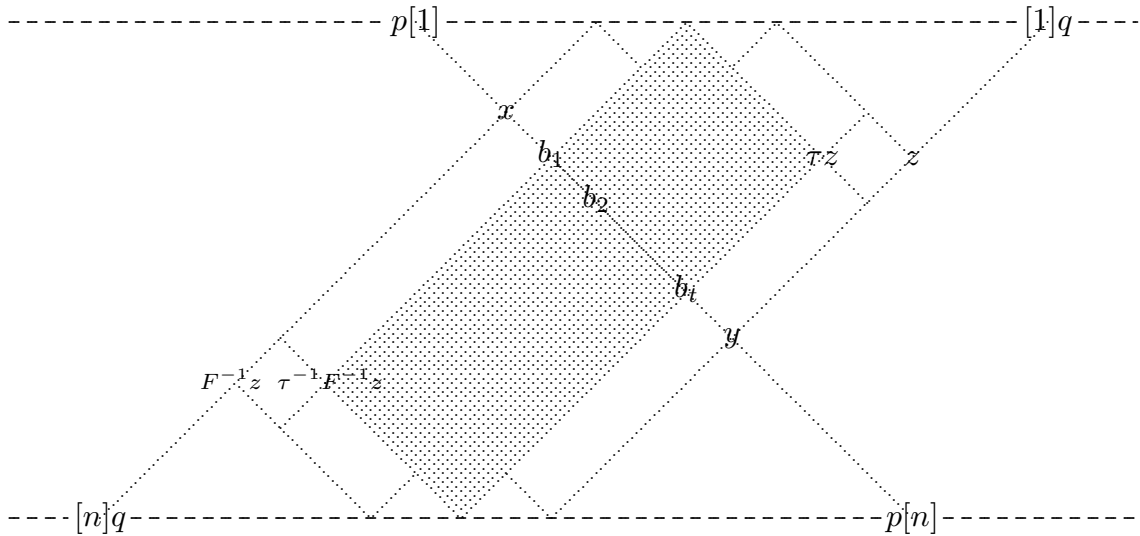
$$

$$
{\beginpicture
\setcoordinatesystem units <.5cm,.5cm>
\setdashes <1mm> 
\plot -2 0  -.7 0 /
\plot .7 0  17.3 0 /
\plot 18.7 0  23  0 /
\plot -2 11  6.3 11 /
\plot 7.7 11  20.3 11 /
\plot 21.7 11  23 11 /

\put{$p[1]$} at 7 11
\put{$[1]q$} at 21 11
\put{$p[n]$} at 18 0
\put{$[n]q$} at 0 0

\setdots <.7mm>
\plot 0 0  11 11  12 10  13 11  14 10  15 11  18  8 /
\plot 21  11  10 0  9 1  8 0  7 1  6  0  3 3 /
\plot 7 11  18 0 /
\plot 4 4  7 1 /
\plot 9 1  17 9 /
\plot 17 7  14 10 /
\plot 4 2    12 10 /

\multiput{} at -1 0  22 11 /

\put{$x$} at 9 9
\put{$y$} at 14 4
\put{$z$} at 18 8
\put{$\tau z$} at 16 8
\put{$\ssize F^{-1}z$} at 2.9 3.1
\put{$\ssize \tau^{-1}F^{-1}z$} at 5.1 3.1
\put{$b_1$} at 10 8.1
\put{$b_2$} at 11 7.1

\put{$b_t$} at 13 5.1

\setshadegrid span <.5mm>
\vshade 5 3 3 <z,z,,> 8 0 6  <z,z,,> 13 5 11  <z,z,,> 16 8 8 /

\endpicture}
$$

According to the wing lemma,
we know that 
$$
 f(z) = -\min(f(b_i)\mid 1\le i \le t) < 0
$$

But starting with $x$ and $y$, we may also look at the wing
with corners 
$$
 [n]q,\ x,\ p[n],
$$ 
and use the dual argument: the dual of the wing lemma
concerns the vertex $F^{-1}z$ (as well as the vertices between
$F^{-1}z$ and $x$), it yields 
$$
 f(F^{-1}z) = -\min(f(b_i)\mid 1\le i \le t).
$$
This shows that 
$$
 f(z) = f(F^{-1}z).
$$
Also, the rectangle rule for $F^{-1}z$ (or the dual 
rectangle rule for $z$) assert that $f$ is bounded from below
by $f(z)^- = -f(z)$ on the rectangle starting
with $\tau^{-1}F^{-1}z$ and ending with $\tau z$ (the shaded
area). 

Using induction as well as duality, we see that 
$f(F^az) = f(z)$ for all $a\in \Bbb Z$. Also, it follows that
$$
 f(z)^-h_z \le f.
$$
	\medskip

Proof of Theorem 4. Choose some slice $\Cal S$. Given a function $g$
on the set of vertices of $\Gamma,$ we write
$$
 |g|_{\Cal S} = \sum_{s\in \Cal S} |g(x)|.
$$
thus $ |g|_{\Cal S} = 0$ if and only if $g(x) = 0$ for all $x\in \Cal S.$
In case $g$ is cluster-additive, we know from section 1 that
$ |g|_{\Cal S} = 0$ if and only if $g$ is the zero function.
	\medskip
We want to show any cluster-additive function $f$ on $\Gamma$
is a non-negative linear combination of cluster-hammock functions. 
We use induction
on $|f|_{\Cal S}.$ If $|f|_{\Cal S} = 0$, then $f$ is the zero function.

Now assume that $|f|_{\Cal S} > 0$. According to the assertion (5) in section
1, there is some vertex $z$ with $f(z) < 0.$ 

According to (2), we know that $h_z \le f(z)^-h_z \le f$.
We see by Theorem 2 that $f-h_z$ is cluster-additive again, and
$|f-h_z|_{\Cal S} < |f|_{\Cal S}.$ Thus, by induction,
$f-h_z$ is a non-negative linear combination of cluster-hammock
functions and then also $f = (f-h_z)+h_z$ is a 
non-negative linear combination of cluster-hammock
functions. This completes the proof.
	\bigskip\bigskip
{\bf 10. Cluster-tilted algebras.}
     \medskip

Let $A$ be a finite-dimensional
hereditary $k$-algebra ($k$ an algebraically closed field). 
Let $T$ be a tilting $A$-module, $\Cal T$ the set of isomorphism
classes of indecomposable direct summands of $T$,
and $F\Cal T$ the union of the $F$-orbits which contain elements of $\Cal T.$
Let $B$ be the opposite endomorphism ring of $T$ in the cluster category 
$\Cal C = D^b(\mod A)/F$ (see [B-T]), thus $B$ is a cluster-tilted
algebra. 

Define a function $d_T$ on the Auslander-Reiten quiver $\Gamma$ of 
$D^b(\mod A)$ as follows: Consider the projection
$$
 D^b(\mod A) \longrightarrow D^b(\mod A)/F = \Cal C_A \longrightarrow \Cal C_A/\langle T\rangle = \mod B,
$$
and denote it by $\pi$.

Let $y$ be a vertex of $\Gamma$, thus the isomorphism class
of an indecomposable object in $D^b(\mod A)$. If $y$ is not in $F\Cal T$, 
then $\pi(y)$ is
the isomorphism class of an indecomposable $B$-module and we denote
by $d_T(y)$ its $k$-dimension. On the other hand, if the $F$-orbit of $y$
contains an element $x$ of $\Cal T$, and $x = [X]$, where $X$ is an indecomposable direct summand of $T$, then let $d_T(x) = n_x$ be the Krull-Remak-Schmidt multiplicity of $X$ in $T$, note that  this is also the $k$-dimension of the corresponding simple $B$-module $S_x$. 
In this way we obtain a function
$$
 d_T\:\Gamma_0 \to \Bbb Z
$$
which obviously is $F$-invariant. 

Of course, instead of looking at the $k$-dimension of the $B$-modules, one may also
consider their length. In this way, one similarly defines the function
$$
 l_T\:\Gamma_0 \to \Bbb Z
$$
with $l_T(y)$ the length of $\pi(y)$ in case $y$ is not in $F\Cal T$,
and with $l_T(y) = -1$ otherwise.
If the tilting module $T$ is multiplicity-free, then $l_T = d_T$. 
For a general tilting module $T$, let $T'$ be multiplicity-free with the same
indecomposable direct summands as $T$, then $l_T = d_{T'}.$

	  \medskip 
{\bf Lemma.} {\it The function $d_T$ on $\Gamma$ is cluster-additive and we have}
$$
 d_T = \sum\nolimits_{x\in \Cal T} n_xh_x.
$$
	\medskip
Proof: Let us consider the mesh of $\Bbb Z\Delta$ ending in $z$,
say with arrows $y_i \to z$, $1\le i \le s$. We assume that 
the vertices $y_{r+1},\dots, y_s$ belong to $F\Cal T$,
and $y_1,\dots, y_r$ not.

First, consider the case that neither $z$ nor $\tau z$ belong to
$F\Cal T,$ thus we may consider the Auslander-Reiten sequence
ending in $Z$. By the assumption on the $y_i$, we see that 
the Auslander-Reiten sequence has the form
$$
 0 \to \tau Z \to \bigoplus_{i=1}^r Y_i^{m_i} \to Z \to 0,
$$
with indecomposable $B$-modules $Z$ and $Y_i$ such that
$[Z] = z,$ $[Y_i] = y_i$ and where $m_i = m_{y_i,z}$. It follows that
$$
\align
 d_T(z) + d_T(\tau z) &= \dim Z + \dim \tau Z \cr
  &= \sum_{i=1}^r m_i\dim Y_i = \sum_{i=1}^r m_{y_i,z} d_T(y_i) \cr
  &=  \sum_{i=1}^s m_{y_i,z} d_T(y_i)^+,
\endalign
$$
since $d_T(y_i) < 0$ for $r+1\le i \le s.$

Next, let $\tau z$ belong to $F\Cal T$, thus $z$ is a projective
vertex, say $z = [Z]$ for some indecomposable projective
$B$-module $Z$. By the assumption on the $y_i$, the radical
of $Z$ has the form $\rad Z = \bigoplus_{i=1}^r Y_i^{m_i}$ and
$Z/\rad Z$ has dimension $n_z.$ This means
$$
\align
 d_T(z) + d_T(\tau z) &= \dim Z - \dim Z/\rad Z  = \dim \rad Z \cr
  &= \sum_{i=1}^r m_i\dim Y_i = \sum_{i=1}^r m_{y_i,z} d_T(y_i) \cr
  &=  \sum_{i=1}^s m_{y_i,z} d_T(y_i)^+.
\endalign
$$

Finally, we have to consider the case where $z$ belongs to
$F\Cal T.$ This case is dual to the previous one, now $\tau z = [X]$
for some indecomposable injective $B$-module $X$ and the socle
of $X$ has dimension $n_z.$ 

The Jordan-H\"older theorem for $\mod B$ shows that 
$d_T$ is just the sum of the various functions $n_xh_x$
with $x\in \Cal T$; namely, if $y$ is a vertex of $\Gamma$, 
such that $\pi(y)$ is
the isomorphism class of an indecomposable $B$-module $N$,
then $h_x(y)$ is just the Jordan-H\"older multiplicity of the simple
$B$-module $S_x$ in $N$.

	\bigskip
In the Dynkin case, we can use the cluster-algebras in order to
prove our conjecture for an $F$-invariant cluster-additive function $f$
provided two conditions on the position of the vertices $x$ with
$f(x) \le 0$ are satisfied. 
	\medskip
{\bf Proposition.} {\it Let $f$ be a cluster-additive function on 
$\Gamma = \Bbb Z\Delta$ with $\Delta$ a Dynkin quiver.
Assume that $f$ is $F$-invariant
and that there is a tilting set $\Cal T$ with the following two properties:
\item{\rm (a)} If $x$ belongs to $\Cal T$, then $f(x) \le 0$.
\item{\rm (b)} If $f(x) < 0$, then $x$ belongs to the $F$-orbit of an
element of $\Cal T$.

\noindent
Then $f$ is a non-negative linear combination of cluster-hammock functions.}
	\medskip
Proof: We identify $\Gamma = \Bbb Z\Delta$ with the Auslander-Reiten quiver of
$D^b(\mod A)$ where $A$ is a finite-dimensional hereditary algebra and where 
$T$ is a tilting $A$-module such that $\Cal T$ is the set of isomorphism
classes of indecomposable direct summands of $T.$ Let $B$
be the opposite endomorphism ring of $T$ in 
$\Cal C_A = D^b(\mod A)/F$.
We form the
factor category $D^b(\mod A)/\langle F^iT\mid i\in \Bbb Z\rangle$,
this is the module category of a Galois cover $\widetilde B$
of $B$ (with Galois group $\Bbb Z$). Thus, the Auslander-Reiten
quiver 
$\Gamma' = \Gamma(\widetilde B)$ of $\mod \widetilde B$ is the translation subquiver 
obtained from $\Gamma$ by deleting the $F$-orbits of the vertices in $\Cal T$.

Denote by $f'$ the restriction of $f$ to $\Gamma'$.
By assumption (b),
$f'$ takes values in $\Bbb N_0$, is cluster-additive, thus additive on $\Gamma(\widetilde B)$ and
$F$-invariant; thus it induces an additive function $f''$ on 
$\Gamma'' = \Gamma(B) = \Gamma(\widetilde B)/F = \Gamma'/F$
with values in $\Bbb N_0.$ 
According to Butler [Bu], $f''$ is additive on all exact sequences,
thus it is a linear combination of the ``hammock functions'' $h''_p$
for $\mod B$,  where $p$ runs through the set of indecomposable
projective $B$-modules. If we compose these functions $h''_p$ 
with the projection $\Gamma' \to \Gamma'/F = \Gamma''$, we obtain
just the restriction of $h'_p$ to $\Gamma'$, where $p = \tau^{-1}x$
for some $x\in \Cal T.$ Thus, there are integers $n_p$ such that
$$ 
 f'' = \sum_{p} n_ph''_p,
$$  
and therefore
$$
 f|\Gamma' =  f' = \sum_p n_ph'_p|\Gamma'.
$$
If $P'$ is an indecomposable projective $B$-module with isomorphism
class $p'$ and $S'$ is its top (a simple $B$-module), then
$$
 n_{p'} = \sum_p n_ph''_p([S']) = f''([S']) \ge 0
$$
(here we use that $f'$ takes values in $\Bbb N_0$), thus all the
coefficients $n_p$ are non-negative. 

We have seen that $f$ and $h = \sum_p n_ph'_p$ coincide on $\Gamma'$,
it remains to show that they also coincide on $\Cal T.$ 
Let $x\in \Cal T$, thus $p = \tau^{-1}x$ is in $\Gamma'_0$ and 
$$
\align
 f(x) = &= -f(p) + \sum_{y\in \Gamma_0} m_{y,p} f(y)^+  \cr
        &= -f(p) + \sum_{y\in \Gamma_0'} m_{y,p(s)}f(y)^+  \cr
        &= -h(p) + \sum_{y\in \Gamma_0'} m_{y,p(s)}h(y)^+  \cr
        &= -h(p) + \sum_{y\in \Gamma_0} m_{y,p(s)}h(y)^+  \cr
        &= h(x),
\endalign
$$
where we have used that both $f$ and $h$ are cluster-additive,
that the coincide on $\Gamma'$ and have positive values only on
vertices in $\Gamma'_0$ (condition (a)). 
This completes the proof that $f = h.$

	\bigskip
If a cluster-additive function on $\Gamma$ 
is a non-negative linear combination of cluster-hammock functions, 
then also the following properties are satisfied:
   \medskip
(d) Always, $f = d_T$ for some partial tilting module $T$.
    \medskip
(e) If $f$ takes values in $\{-1\}\cup \Bbb N_0$, then $f= d_T$ for some multiplicity free partial tilting module $T$, if $f$ 
takes values in $\Bbb Z\setminus \{0\}$, then $f= d_T$ for some tilting module $T$.
	\bigskip
We end this section by giving an interpretation of the 
exchange property of cluster-tilting objects in a cluster category in terms of the
cluster-hammock functions. 
Thus, suppose that we deal with a tilting set $\Cal T$ in $\Bbb Z\Delta$, where
$\Delta$ is a Dynkin quiver. 
Let us look at the hammock $h_x$  for some $x\in \Cal T.$ Let $\Cal T' = 
\Cal T \setminus \{x\}$. We claim that there are precisely two $F$-orbits of 
vertices of $\Gamma$ which are not in the support of any function $h_{y}$ with
$y\in \Cal T'$. Of course, one of these vertices is $x$ itself, since $h_y(x) = 0$
for all $y\in \Cal T'.$ In order to find the other orbit, we only have
to consider the vertices $z$ which do not belong to $F\Cal T$. As above, we know
that  $\pi(z)$ is the isomorphism class $[N]$ of an indecomposable $B$-module, say $N$.
Now either $[N] = [S_x]$, then indeed $h_y(z) = 0$ for all $y\in \Cal T'$ (since
$N$ has no composition factor of the form $S_y$),
or else $N$ is not isomorphic to $S_x$, but then $N$ has at least one composition
different from $S_x$, say $S_y$ with $y\in \Cal T',$
and therefore $h_y(z) \neq 0.$ This shows that the second orbit consists of the
vertices $z$ such that $\pi(z) = [S_x]$.  
(But a warning is necessary: the position of $z$ 
with $\pi(z) = [S_x]$ in the support
of $h_x$ does not only depend on $h_x$ itself, as already the case $\Bbb A_2$ shows.)
	\bigskip\bigskip
{\bf 11. Final Remarks.}
     \medskip
{\bf 1.} The main results and conjectures of this note concern the translation quivers $\Bbb Z\Delta$ with
$\Delta$ a simply laced Dynkin diagram. But there is no problem to extend  the considerations to
the case of an arbitrary (not necessarily simply laced) 
Dynkin diagram. In order to do so, we need the
notion of a valued translation quiver.

A {\it valued translation quiver} $\Gamma = (\Gamma_0,\Gamma_1,\tau,v)$ is
given by a translation quiver $(\Gamma_0,\Gamma_1,\tau)$ with the property
that there is at most one arrow $x\to y$ for any pair $x,y$ of vertices
and a function $v\:\Gamma_1 \to \Bbb N_1$ such that $v_{\tau x,\tau y} =
v_{x,y}$ for all arrows $x\to y$ (we write $v_{x,y}$ or just $v_{xy}$
instead of $v(x\to y)$).
In case
$v_{\tau z,y} = v_{y,z}$ for all arrows $y\to z$ with $z$ not projective.
then we may consider $(\Gamma_0,\Gamma_1,\tau,v)$ as an ordinary
translation quiver by replacing any arrow $x\to y$ by $v_{x,y}$ arrows.

For example, the valued translation quiver $\Bbb Z\Bbb B_3$ has the following
form (in such pictures it is sufficient to add the number $v_{xy}$ to an arrow $x\to y$ only in case $v_{xy} \ge 2$):
$$
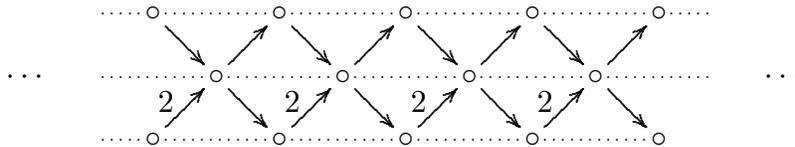
{\beginpicture
\setcoordinatesystem units <.7cm,.7cm>
\multiput{} at 0 -.3  8 2.3 /
\multiput{$\circ$} at 0 0  2 0  4 0  6 0  8 0  
                      0 2  2 2  4 2  6 2  8 2
                      1 1  3 1  5 1  7 1 /
\arr{0.2 1.8}{0.8 1.2}
\arr{0.2 0.2}{0.8 0.8}
\arr{1.2 1.2}{1.8 1.8}
\arr{1.2 0.8}{1.8 0.2}
\put{$2$} at 0.2 0.6

\arr{2.2 1.8}{2.8 1.2}
\arr{2.2 0.2}{2.8 0.8}
\arr{3.2 1.2}{3.8 1.8}
\arr{3.2 0.8}{3.8 0.2}
\put{$2$} at 2.2 0.6

\arr{4.2 1.8}{4.8 1.2}
\arr{4.2 0.2}{4.8 0.8}
\arr{5.2 1.2}{5.8 1.8}
\arr{5.2 0.8}{5.8 0.2}
\put{$2$} at 4.2 0.6

\arr{6.2 1.8}{6.8 1.2}
\arr{6.2 0.2}{6.8 0.8}
\arr{7.2 1.2}{7.8 1.8}
\arr{7.2 0.8}{7.8 0.2}
\put{$2$} at 6.2 0.6
\setdots <1mm>
\plot -.8 2  -.2 2 /
\plot 0.2 2  1.8 2 /
\plot 2.2 2  3.8 2 /
\plot 4.2 2  5.8 2 /
\plot 6.2 2  7.8 2 /
\plot 8.2 2  8.8 2 /

\plot -.8 0  -.2 0 /
\plot 0.2 0  1.8 0 /
\plot 2.2 0  3.8 0 /
\plot 4.2 0  5.8 0 /
\plot 6.2 0  7.8 0 /

\plot -.8 1  .8 1 /
\plot 1.2 1  2.8 1 /
\plot 3.2 1  4.8 1 /
\plot 5.2 1  6.8 1 /
\plot 7.2 1  8.8 1 /

\multiput{$\cdots$} at -2 1  10 1 /
\endpicture}
$$

The valued translation quiver $\Gamma = (\Gamma_0,\Gamma_1,\tau,v)$
is said to be {\it stable,} if $(\Gamma_0,\Gamma_1,\tau)$ is stable.
Given a stable valued translation quiver 
$\Gamma$, a function $f\:\Gamma_0 \to \Bbb Z$ should be called 
{\it cluster-additive} provided
$$
 f(z)+f(\tau z) = \sum\nolimits_{y\in \Gamma_0} v_{yz}f(y)^+, \quad \text{for all\ } z\in \Gamma_0.
$$

	\bigskip
{\bf 2.} We should stress that cluster-additive functions are 
definitely also of interest when dealing with stable
translation quivers which are not related to  translation quivers
of the form $\Bbb Z\Delta$ with  $\Delta$ a finite directed quiver.
Examples of  cluster-additive functions
on the translation quiver $\Bbb Z\Bbb D_\infty$ (as well as on
$\Bbb Z\Bbb A_\infty^\infty$)  have been exhibited in [R].

   \bigskip\bigskip
\vfill\eject
{\bf References.}

\frenchspacing
     \medskip
\item{[BGP]} I. N. Bernstein, I. M. Gelfand and V. A. Ponomarev,
Coxeter functors and Gabriel's theorem.
Uspechi Mat. Nauk. 28 (1973), 19-33; Russian Math. Surveys 29 (1973), 17-32.

\item{[Br]} S. Brenner: A combinatorial characterisation of finite Auslander-Reiten quivers. In: Representation Theory I. Finite
Dimensional Algebras (ed. V. Dlab, P. Gabriel, G. Michler).
Springer LNM 1177 (1986), 13-49.

\item{[B-T]} A. B. Buan, R. J. Marsh, M. Reineke, I. Reiten, G. Todorov, Tilting theory and cluster
combinatorics, Adv. Math. 204 (2006), 572.618.

\item{[BMR]} A. B. Buan, R. J. Marsh, I. Reiten, Cluster-tilted algebras, Trans. Amer. Math. Soc. 359
(2007), 323.332.

\item{[Bu]} M.C.R. Butler, Grothendieck groups and almost split sequences.
 In: Integral Representations and Applications (ed. K.W. Roggenkamp).
Springer LNM 882 (1981), 357-368.

\item{[G]} P. Gabriel. Auslander-Reiten sequences and representation-finite
algebras. In: Representation Theory I, Springer LNM 831 (1980), 1-71.

\item{[H]} D. Happel: Triangulated categories in the representation
theory of finite-dimensional algebras. LMS Lecture Note Series 119.
Cambridge (1988).

\item{[HPR]} D. Happel, U. Preiser, C.M. Ringel:
Vinberg's characterization of Dynkin diagrams using subadditive functions with application to DTr-periodic modules. In: Representation Theory II, 
Springer LNM 832 (1980), 280-294.

\item{[R]} C.M. Ringel, The minimal representation-infinite algebras which are special biserial. To appear in: Representations of Algebras and Related Topics,  European Math. Soc. Publ. House, Z\"urich
(ed. A. Skowro\'nski, K. Yamagata) (2011)

\item{[RV]} C.M. Ringel, D. Vossieck: 
  Hammocks. Proc. London Math. Soc. (3) 54 (1987), 216-246.

\item{[S]} Scheuer: More hammocks. In: Topics in Algebra. Banach Center Publications
26, Part I, 493-498

\bigskip\bigskip

{\rmk Fakult\"at f\"ur Mathematik, Universit\"at Bielefeld \par
POBox 100\,131, \ D-33\,501 Bielefeld, Germany \par
e-mail: \ttk ringel\@math.uni-bielefeld.de \par}

\bye